\documentclass[12pt]{amsart}
\usepackage{amsmath,amssymb,latexsym}
\usepackage{nicefrac}
\usepackage{pgfplots}
\usepackage{subcaption}
\usepackage{hyperref}
\theoremstyle{plain}             
\newtheorem{theorem}{Theorem}[section]

\newtheorem{definition}[theorem]{Definition}

\newcommand{\isdef}{\mathrel{\mathrel{\mathop:}=}}

\newcommand{\spn}{\operatorname{span}}
\newcommand{\bs}{\boldsymbol}
\renewcommand{\d}{\operatorname{d}\!}
\newcommand{\dn}{\operatorname{dn}}
\newcommand{\cn}{\operatorname{cn}}
\newcommand{\sn}{\operatorname{sn}}
\definecolor{green}{RGB}{28,172,0}
\definecolor{lila}{RGB}{170,55,241}

\begin{document}
\title[Simulation of random fields on surfaces]{Simulation of Gaussian 
random fields on surfaces using the isogeometric finite element method}
\author{Helmut Harbrecht \and Florian Sonderegger \and Remo von Rickenbach}
\address{Helmut Harbrecht, Florian Sonderegger, Remo von Rickenbach,
Department of Mathematics and Computer Science,
University of Basel, Spiegelgasse 1, 4051 Basel, Switzerland}
\email{\{helmut.harbrecht,f.sonderegger,remo.vonrickenbach\}@unibas.ch}

\maketitle
\begin{abstract}
We are concerned with the fast simulation of random fields 
on closed surfaces in $\mathbb{R}^3$ which are generated 
by the (Whittle-) Mat\'ern class of covariance functions. To 
this end, we solve the underlying fractional stochastic partial 
differential equation with additive white noise by using an 
isogeometric finite element method on the surface in combination 
with the Balakrishnan integral representation of the solution. 
The solution of the underlying linear system of equations is 
performed by means of a geometric multigrid method that 
naturally underlies the isogeometric approach. Numerical 
results are presented to demonstrate the approach.
\end{abstract}

\section{Introduction}
A commonly used category of statistical models is that 
of Gaussian processes, characterized by their mean and 
covariance structure. This involves the generation of correlated 
random data that follow a Gaussian distribution. In this article, 
we are interested in simulating such \emph{Gaussian random 
fields} on surfaces of three-dimensional domains. They have 
been utilized in a variety of fields, for example in engineering, 
geostatistics, cosmological data analysis, and biomedical 
image analysis, compare \cite{application1,application2,
application3} and the references therein.

To generate realizations of a specific Gaussian random field,
we pursue with an approach that has its origin in computational 
spatial statistics \cite{HHKS,Chris,Annika3,Rue1,Rue2}. This approach 
allows the efficient sampling for Gaussian random fields arising from 
the (Whittle-) Mat\'ern class of covariance functions. As observed for 
$\mathbb{R}^d$ in \cite{Whittle1,Whittle2}, these covariance 
functions are given by
\begin{equation}\label{eq:kernel}
  \mathcal{K}_{\nu}({\bs x},{\bs y})\isdef
  \frac{\sigma^2}{2^{\nu-1}\it{\Gamma}(\nu)} r^{\nu} K_{\nu} (r),
  \quad \text{where}\ r \isdef \kappa\|{\bs x}-{\bs y}\|,
\end{equation}
where $\nu$ is the smoothness parameter, determining 
the spatial (Sobolev) smoothness of the Mat\'ern field, 
$\kappa>0$ relates to the correlation length, $\it{\Gamma}$ 
is the Riemannian gamma function and $K_\nu$ is the 
modified Bessel function of the second kind, see \cite{MAT} 
for example. The marginal variance $\sigma^2$ in 
\eqref{eq:kernel} is given by 
\begin{equation}\label{eq:sigma}
\sigma^2 = \frac{\it{\Gamma}(\nu)}{\it{\Gamma}(\nu+\nicefrac{d}{2})}
(4\pi)^{-\nicefrac{d}{2}}\kappa^{-2\nu}.
\end{equation}

In the situation of a two-dimensional closed surface
$\Gamma := \partial D$ of a given domain $D\subset
\mathbb{R}^3$, we are in the situation of $d=2$ and 
the distance $\|{\bs x}-{\bs y}\|$ in \eqref{eq:kernel} 
should be the geodesic distance which is not easy to 
calculate on general surfaces. Therefore, we use the 
characterization of the random field as the solution of 
the following fractional partial differential equation with 
additive white noise:
\begin{equation}\label{eq:PDE}
(\kappa^2 I+\Delta_\Gamma)^\beta u = f\ \text{on $\Gamma$}.
\end{equation}
Herein, $\Delta_\Gamma$ denotes the Laplace-Beltrami 
operator on the given surface $\Gamma\subset\mathbb{R}^3$
while the fractional exponent $\beta$ is related to the smoothness 
parameter $\nu$ in \eqref{eq:kernel} via $\nu = 2\beta-1$. 

For $\beta\in(0,1)$, the solution of the fractional partial differential 
equation can be computed by using the integral identity
\begin{equation}\label{eq:integral}
(\kappa^2 I-\Delta_\Gamma)^\beta u 
= \frac{\sin(\pi\beta)}{\pi}\int_0^\infty t^{-\beta} (t+\kappa^2 -\Delta_\Gamma)^{-1} u\d t,
\end{equation}
which has been proven in \cite{Balakrishnan}. In particular,
this improper integral can be efficiently approximated by a sinc 
quadrature, see \cite{Lund} for details. For $\beta > 1$, one 
splits $\beta = n+\beta'$ with the integer part $n\in\mathbb{N}$ 
being successively computable by solving for the Laplace-Beltrami 
operator and the fractional part $\beta'\in(0,1)$ being computable 
by means of \eqref{eq:integral}, compare \cite{Chris}. However,
for $\beta'$ being close to $0$ or $1$, the convergence of the sinc
quadrature becomes quite slow and we propose an improved
splitting in Section~\ref{sct:fractional}.

The evaluation of the integrand in \eqref{eq:integral} 
at a sinc quadrature point amounts to solving a partial 
differential equation for the Laplace-Beltrami operator.
In contrast to \cite{Bonito,Annika2,Annika1}, we apply here
the methodology from \emph{isogeometric analysis} (IGA) 
to discretize this partial differential equation. This means 
that we employ a high-order or even exact representation of 
the underlying surface in terms of \emph{non-uniform rational 
B-splines} (NURBS) in combination with a high-order spline 
based discretization of the partial differential equation. This
enables to exploit the higher order smoothness of the underlying
Gaussian random field when $\beta$ is getting larger.

IGA has been introduced in \cite{HCB05} in order to incorporate 
simulation techniques into the design workflow of industrial development
and thus allows to deal with surfaces in a straightforward manner. In 
addition, the naturally emerging sequence of nested approximation 
spaces can directly be employed in a multigrid method. Therefore, 
the isogeometric finite element method is the method of choice for 
the problem at hand. In particular, we can rely on the \verb|C++|-library 
\verb+bembel+ which provides an implementation of the isogeometric 
boundary element method, naturally including isogeometric 
finite elements on surfaces, see \cite{bembel,DHK+20,DHK+18}
and also \cite{huang_isogeometric_2022}.

The rest of this article is organized as follows. Section~\ref{sec:IGA}
recapitulates the basic concepts from isogeometric analysis and 
introduces the discretization spaces that will be used later on. 
In Section~\ref{sct:discretization}, we introduce the discretization 
of the fractional boundary value problem under consideration and
develop related solution algorithms. In Section~\ref{sec:numerix},
we present numerical results to validate the applicability and 
feasibility of the proposed approach for the simulation of Gaussian 
random fields arising from (Whittle-) Mat\'ern covariance functions. 
Finally, in Section~\ref{sct:conclusio}, we draw the article's 
conclusion.

\section{Isogeometric analysis}\label{sec:IGA}
\subsection{B-splines}
We shall give a brief introduction to the basic concepts 
of isogeometric analysis, starting with the definition of the 
B-spline basis, followed by the description of the surface
under consideration by using NURBS. The original definitions 
(or equivalent notions) and proofs, as well as basic algorithms, can 
be found in most of the standard spline and isogeometric literature
\cite{Cottrell_2009aa,HCB05,Piegl_1997aa,Schumaker_1981aa,Lee_1996aa}.

\begin{figure*}
\begin{subfigure}{.5\textwidth}
\begin{tikzpicture}
\begin{axis}[
xmin = -.1,
xmax = 1.1,
ymin = -.1,
ymax = 1.1,
width=1\columnwidth,
height=.66\columnwidth,
grid=major,
legend style={
at={(.5,1)},
anchor=south}
]
\addplot[red,ultra thick,mark=none,domain = 0:1/3] {1};
\addplot[teal,ultra thick,mark=none,domain = 1/3:2/3] {1};
\addplot[blue,ultra thick,mark=none,domain = 2/3:1] {1};
\end{axis}
\end{tikzpicture}
\caption{$p=0$, $\Xi=[0,1/3,2/3,1]$.}
\end{subfigure}\!\!
\begin{subfigure}{.5\textwidth}
\begin{tikzpicture}
\begin{axis}[
xmin = -.1,
xmax = 1.1,
ymin = -.1,
ymax = 1.1,
width=1\columnwidth,
height=.66\columnwidth,
grid=major,
legend style={
at={(.5,1)},
anchor=south}
]
\addplot[red,ultra thick,mark=none,domain = 0:1/3] {3*-(x-1/3)};
\addplot[teal,ultra thick,mark=none,domain = 0:1/3] {3*(x)};
\addplot[blue,ultra thick,mark=none,domain = 1/3:2/3] {3*(x-1/3)};
\addplot[orange,ultra thick,mark=none,domain = 2/3:1]{3*(x-2/3)};
\addplot[teal,ultra thick,mark=none,domain = 1/3:2/3] {(1-3*x)+1};
\addplot[blue,ultra thick,mark=none,domain = 2/3:1] {3*(1-x-1/3)+1};
\end{axis}
\end{tikzpicture}
\caption{$p=1$, $\Xi=[0,0,1/3,2/3,1,1]$.}
\end{subfigure}\\[.5cm]
\begin{subfigure}{\textwidth}
\begin{tikzpicture}
\begin{axis}[
xmin = -.1,
xmax = 1.1,
ymin = -.1,
ymax = 1.1,
width=1\columnwidth,
height=.3397\columnwidth,
grid=major,
legend style={
at={(1,.5)},
anchor=west}
]
\addplot[red,ultra thick,mark=none,domain = 0:1/3]{(3*(x)-1)^2} ;
\addplot[teal,ultra thick,mark=none,domain = 0:1/3]{2*(3*(x))*(1-3*x)+.5*(3*x)^2};
\addplot[blue,ultra thick,mark=none,domain = 0:1/3]{.5*(x*3)^2};
\addplot[orange,ultra thick,mark=none,domain = 2/3:1]{2*(3*((x-2/3)))*(1-3*(x-2/3))+.5*(1-3*(x-2/3))^2};
\addplot[brown,ultra thick,mark=none,domain = 2/3:1]{(3*(x-2/3))^2};
\addplot[teal,ultra thick,mark=none,domain = 1/3:2/3]{.5*(1-3*(x-1/3))^2 };
\addplot[blue,ultra thick,mark=none,domain = 1/3:2/3]{-(((x*3)-1)-1)*((x*3)-1)+.5};
\addplot[blue,ultra thick,mark=none,domain = 2/3:1]{.5*(3-(x*3))^2};
\addplot[orange,ultra thick,mark=none,domain = 1/3:2/3]{.5*(3*(x-1/3))^2 };
\end{axis}
\end{tikzpicture}
\caption{$p=2$, $\Xi=[0,0,0,1/3,2/3,1,1,1]$.}
\end{subfigure}
\caption{B-spline bases for $p=0,1,2$ and open knot 
vectors with interior knots $1/3$ and $2/3$.}\label{fig:splines}
\end{figure*}
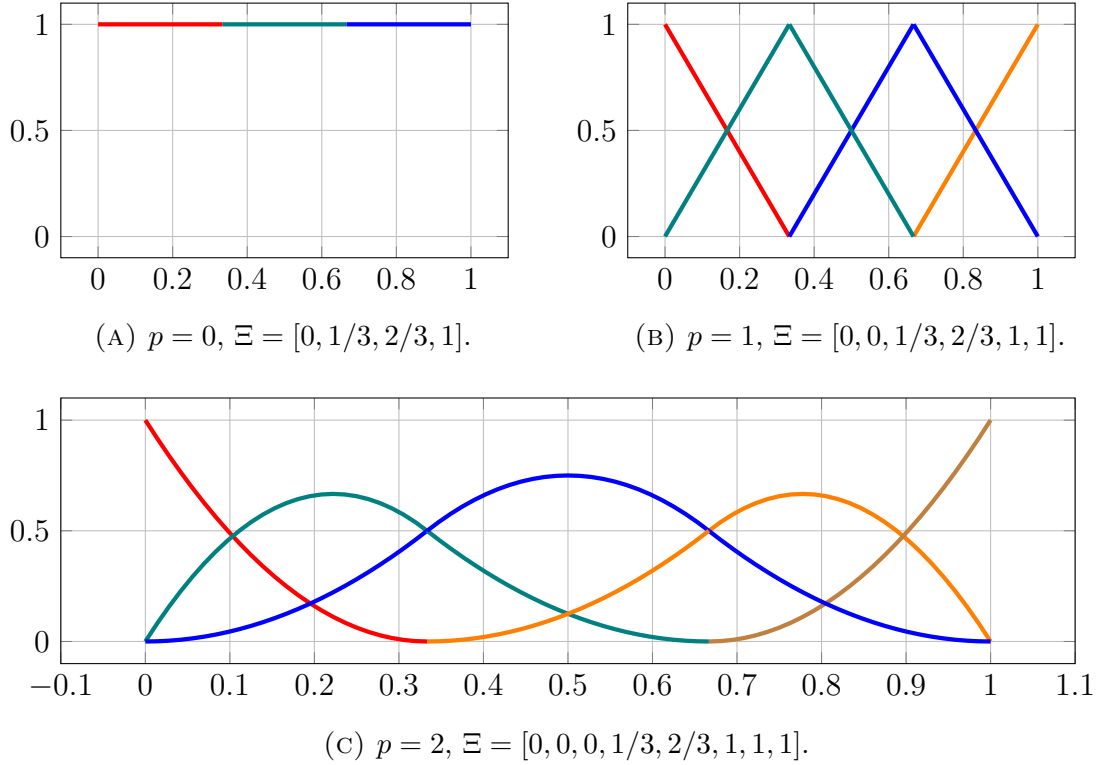

\begin{definition}
Let $0\leq p\leq k$. We define a \emph{$p$-open knot 
vector} as a set  $\Xi := \{\xi_0,\dots, \xi_{n+p}\}$ such that
\[
\underbrace{\xi_0 = \cdots =\xi_{p}}_{=0}<\xi_{p+1}<\cdots <\xi_{k-1}<
	\underbrace{\xi_{k}=\cdots =\xi_{k+p}}_{=1},
\]
where $k$ denotes the number of control points. The 
associated basis functions are given by $\{b_\ell^p\}_{\ell=0}^{k-1}$ 
for $p=0$ as
\[
  b_\ell^0(x) =\begin{cases}
	1, & \text{if }\xi_\ell\leq x<\xi_{\ell+1}, \\
	0, & \text{otherwise},\end{cases}
\]
and for $p>0$ via the recursive relationship
\[
  b_\ell^p(x) = \frac{x-\xi_\ell}{\xi_{\ell+p}-\xi_j}b_\ell^{p-1}(x) 
  +\frac{\xi_{\ell+p+1}-x}{\xi_{\ell+p+1}-\xi_{\ell+1}}b_{\ell+1}^{p-1}(x),
\]
see~Figure \ref{fig:splines}. A \emph{spline} is then 
defined as a function
\[
f(x) = \sum_{0\le\ell< k}p_\ell b_\ell^p(x),
\]
where $\{p_\ell\}_{\ell=0}^{k-1}\subset\mathbb{R}$ 
denotes the set of \emph{control points}. If one sets 
$\{{\bs p}_\ell\}_{\ell=0}^{k-1}\subset\mathbb{R}^d$
with $d\ge 2$, then $f$ will be called a \emph{spline curve}.
\end{definition}

Having the spline functions at hand, we can introduce
the spline spaces which serve as fundament for the
definition of the ansatz and test spaces of the finite
element method on the surface. 

\begin{definition}
Let $\Xi$ be a $p$-open knot vector containing $k+p+1$ 
elements. We define the \emph{spline space} $S_{p}(\Xi)$ 
as the space spanned by $\{b_\ell^p\}_{\ell=0}^{k-1}$.
\end{definition}

Finally, we should consider the relationship between the 
spline spaces and the underlying mesh relative to a 
specific mesh size.

\begin{definition}
For a knot vector $\Xi,$ we define the \emph{mesh size} 
$h$ to be the maximal distance 
\begin{equation}\label{eq:h}
h\isdef \max_{\ell=p}^k h_\ell,\quad\text{where}\quad
	h_\ell\isdef\xi_{\ell+1}-\xi_\ell,
\end{equation}
between neighbouring knots.
We call the knot vector \emph{quasi uniform}, when there 
exists a constant $\theta> 0$ such that for all $\ell$ the 
ratio $\nicefrac{h_\ell}{h}$ satisfies $\nicefrac{h_\ell}{h}
\ge \theta$. 
\end{definition}

B-splines on two-dimensional domains are constructed 
by means of the tensor product
\begin{equation}\label{def::tpspline}
  f(x,y) =\sum_{\ell_1=0}^{k_1-1}\sum_{\ell_2=0}^{k_2-1}
  	{\bs p}_{\ell_1,\ell_2} b_{\ell_1}^{p_1}(x) b_{\ell_2}^{p_2}(y)
\end{equation}
with control points ${\bs p}_{\ell_1,\ell_2}\in\mathbb{R}^d$. This 
allows the definition \emph{tensor product spline spaces} 
in accordance with
\[
S_{p_1,p_2}(\Xi_1,\Xi_2).
\]

Throughout this article, we will reserve the letter $h$ for 
the mesh size \eqref{eq:h}. All knot vectors will be assumed to 
be quasi uniform, such that the usual spline theory is applicable,
see~\cite{BuffaActa,Piegl_1997aa,Schumaker_1981aa} 
for example.

\subsection{Isogeometric representation of the surface}
We assume that the boundary $\Gamma$ of the domain
under consideration is closed and Lipschitz continuous. For 
the remainder of this article, we assume that it is given 
patchwise as $\Gamma=\bigcup_{m=1}^n\Gamma_m$, i.e.\ 
that it is induced by smooth diffeomorphisms
\begin{equation}\label{def::geom}
  {\bs F}_m\colon \square = [0,1]^2 \to \Gamma_m \subset \mathbb{R}^3.
\end{equation}
In the spirit of isogeometric analysis, these mappings 
are given by NURBS mappings, i.e.~by
\[
  {\bs F}_m(x,y)\isdef \sum_{\ell_1=0}^{k_1-1}\sum_{\ell_2=0}^{k_2-1}
  \frac{c_{\ell_1,\ell_2} b_{\ell_1}^{p_1}(x) b_{\ell_2}^{p_2}(y) w_{\ell_1,\ell_2}}
  {\sum_{i_1=0}^{k_1-1}\sum_{i_2=0}^{k_2-1} b_{i_1}^{p_1}(x) b_{i_2}^{p_2}(y)}w_{i_1,i_2}
\]
with control points $c_{\ell_1,\ell_2}\in \mathbb{R}^3$ and weights 
$w_{i_1,i_2}>0$. We will moreover require that, for any interface 
$D = \Gamma_{m_1}\cap \Gamma_{m_2} \neq \emptyset$, the 
NURBS mappings coincide, i.e.~that, up to rotation of the reference 
domain, one finds ${\bs F}_{m_1}(\cdot,1) \equiv {\bs F}_{m_2}(\cdot,0)$.

\subsection{Finite element spaces}
The mappings of \eqref{def::geom} give rise to the transformations
\[
  \iota_m(f) \isdef f\circ {\bs F}_m,
\]
which are employed to define discrete spaces patchwise, 
by mapping the space of tensor product B-splines as in 
\eqref{def::tpspline} with
\[
  \Xi_{j,p} \isdef \big\{ \underbrace{0,\dots,0}_{p+1\text{ times}},
  	2^{-j},\dots,  1-2^{-j},\underbrace{1,\dots,1}_{p+1\text{ times}}\big\}
\]
to the geometry. Here, the variable $j$ denotes the level 
of uniform refinement.

Let us first introduce the space of patchwise continuous 
B-splines on $\Gamma$ by
\[
  \mathbb S_{j,p}^2(\Gamma) \isdef \Big\{
  	f\in L^2(\Gamma)\colon f_{|\Gamma_m} \equiv \iota_m^{-1}(g)
	\text{ for some }g\in S_{p,p}(\Xi_{j,p},\Xi_{j,p})\Big\},
\]
It is of dimension $n\cdot(2^j + p )^2$, where $n$ 
denotes the number of patches involved in the description 
of the geometry. For our purpose of discretizing the 
Laplace-Beltrami operator, we however require globally 
continuous B-splines and hence the subspace
\begin{equation}\label{def::space2}
\begin{aligned}
&\mathbb S_{j,p}^0(\Gamma) \isdef \Big\{
f\in C(\Gamma)\colon f_{|\Gamma_m} \equiv \iota_m^{-1}(g)\\
&\hspace*{2.5cm}\text{ for some }g\in S_{p,p}(\Xi_{j,p},\Xi_{j,p})\Big\}
\subset \mathbb S_{j,p}^2(\Gamma).
\end{aligned}
\end{equation}
This space will serve as ansatz and test space in the
Galerkin method to solve the underlying boundary value 
problem. It especially provides the approximation property
\begin{equation}\label{eq:approximation}
  \inf_{u_j\in S_{j,p}^0(\Gamma)}\|u-u_j\|_{L^2(\Gamma)} 
  	\lesssim 2^{-j(p+1)} \|u\|_{H^{p+1}(\Gamma)}
\end{equation}
provided that there holds $u\in H^{p+1}(\Gamma)$.

\section{Discretization}
\label{sct:discretization}
\subsection{Galerkin method}
The evaluation of the integral \eqref{eq:integral} relies 
on the numerical solution of the partial differential equation
\begin{equation}\label{eq:simple solve}
(t+\Delta_\Gamma) u = f\ \text{on $\Gamma$}
\end{equation}
for given $t\in\mathbb{C}$ and some right-hand side $f$. 
To derive the variational formulation of this problem, we let 
\[
\langle u,v\rangle_{\Gamma} = \int_\Gamma 
u({\bs x}) v({\bs x})\d\sigma
\]
denote the surface integral over $\Gamma$.
Then, we the variational formulation reads:
\begin{equation}\label{eq:var_form}
\begin{aligned}
&\text{Find}~u\in H^1(\Gamma)~\text{such that}\\
&\hspace*{1cm}
t \langle u,v\rangle_{\Gamma}
+\langle\nabla_\Gamma u,\nabla_\Gamma v\rangle_{\Gamma}
=\langle f,v\rangle_{\Gamma}\quad\text{for all}~v\in H^1(\Gamma).
\end{aligned}
\end{equation}
Here, $\nabla_\Gamma u$ denotes the surface gradient 
of the function $u$. Given a function $w:\mathbb{R}^3
\to\mathbb{R}$, it is defined in accordance with
\[
\nabla_\Gamma w = \nabla w - \frac{\partial w}{\partial {\bs n}} {\bs n}.
\]

For the discretization of the variational formulation,
we replace the energy space $H^1(\Gamma)$ by
the spline space $\mathbb S_{j,p}^0(\Gamma)$. Thus, the 
discrete variational formulation for \eqref{eq:var_form} reads
\begin{align*}
&\text{Find}~u_h\in \mathbb S_{j,p}^0(\Gamma)~\text{such that}\\
&\hspace*{1cm}
t\langle u_h,v_h\rangle_{\Gamma}
+\langle\nabla_\Gamma u_h,\nabla_\Gamma v_h\rangle_{\Gamma}
=\langle f,v_h\rangle_{\Gamma}\quad\text{for all}~v_h\in \mathbb S_{j,p}^0(\Gamma),
\end{align*}
where $u_h\approx u$ is the Galerkin approximation of the
solution $u$ to \eqref{eq:simple solve}. Setting
$N:=\dim\big(S_{j,p}^0(\Gamma)\big)$ and choosing a basis 
\[
\mathbb S_{j,p}^0(\Gamma)=\spn\{\varphi_1,\ldots,\varphi_N\}
\] 
of the underlying ansatz space leads to the system of 
linear equations
\begin{equation}\label{eq:FEM}
(\kappa^2 {\bs M}+{\bs S}){\bs u}={\bs f}
\end{equation}
with the mass matrix ${\bs M}$, stiffness matrix ${\bs S}$, and right-hand 
side ${\bs f}$ given by
\begin{align*}
{\bs M}=\big[\langle\varphi_j,\varphi_i\rangle_{\Gamma}\big]_{i,j=1}^N,\quad
{\bs S}=\big[\langle\nabla_\Gamma\varphi_j,\nabla_\Gamma\varphi_i\rangle_{\Gamma}\big]_{i,j=1}^N,\quad
{\bs f}=\big[\langle f,\varphi_i\rangle_{\Gamma}\big]_{i=0}^N,
\end{align*}
and ${\bs u}$ being the coefficient vector of $u_h$.

\subsection{Reformulation on the reference domain}\label{sec:bilinearonreference}
Due to the isogeometric representations of the geometry, the 
bilinear forms for the computation of the matrix entries can entirely 
be pulled back to the reference domain \cite{HP13}. We define
the \emph{first fundamental tensor of differential geometry} 
of a mapping ${\bs F}_m$ for $\hat{\bs x} = (x,y)\in [0,1]^2$ by
\[
  {\bs K}_m(\hat{\bs x})\isdef\begin{bmatrix}
  \partial_{x}{\bs F}_m(\hat{\bs x})^\top \partial_{x}{\bs F}_m(\hat{\bs x}) &
  \partial_{x}{\bs F}_m(\hat{\bs x})^\top \partial_{y}{\bs F}_m(\hat{\bs x}) \\
  \partial_{y}{\bs F}_m(\hat{\bs x})^\top \partial_{x}{\bs F}_m(\hat{\bs x}) &
  \partial_{y}{\bs F}_m(\hat{\bs x})^\top \partial_{y}{\bs F}_m(\hat{\bs x}) 
  \end{bmatrix}.
\]
Then, the associated \emph{surface measure} can be 
expressed by
\[
  a_m (\hat{\bs x})\isdef\big\|\partial_{x}{\bs F}_m(\hat{\bs x})
  	\times \partial_{y}{\bs F}_m(\hat{\bs x})\big\|_2
	= \sqrt{\det {\bs K}_m(\hat{\bs x})}.
\]
Thus, the bilinear form related to the mass matrix 
can be recast as
\begin{equation}\label{eq:mass}
\begin{aligned}
  \langle u,v\rangle_{\Gamma} &= \sum_{m=1}^n\int_{[0,1]^2}
  	u\big({\bs F}_m(\hat{\bs x})\big) v\big({\bs F}_m(\hat{\bs x})\big) 
		a_m(\hat{\bs x})\d\hat{\bs x}\\
  &= \sum_{m=1}^n\int_{[0,1]^2} u_m(\hat{\bs x}) v_m(\hat{\bs x}) 
  		a_m(\hat{\bs x})\d\hat{\bs x}	
\end{aligned}
\end{equation}
with the pull-back
\[
u_m(\hat{\bs x})=\iota_m(u)(\hat{\bs x}),\quad
v_m(\hat{\bs y})=\iota_m(v)(\hat{\bs y})
\]
of the ansatz and test functions. 


In order to compute the stiffness matrix, we note 
that the surface gradient of a function $u$ defined 
on the patch $\Gamma_m$ is given by
\[
  \nabla_\Gamma u\big({\bs F}_m(\hat{\bs x})\big) = 
  \begin{bmatrix}\partial_{x}{\bs F}_m(\hat{\bs x}),
  \partial_{y}{\bs F}_m(\hat{\bs x})\end{bmatrix}
  {\bs K}_m^{-1}(\hat{\bs x})
  \begin{bmatrix}
  \partial_{x} u_m(\hat{\bs x})\\
  \partial_{y} u_m(\hat{\bs x})
  \end{bmatrix}.
\]
Hence, we arrive at the identity 
\begin{equation}\label{eq:stiffness}
  \langle \nabla_\Gamma u,\nabla_\Gamma v\rangle_{\Gamma} 
   	= \sum_{m=1}^n\int_{[0,1]^2} 
	\begin{bmatrix}
  	\partial_{x} u_m(\hat{\bs x}),
  	\partial_{y} u_m(\hat{\bs x})
  	\end{bmatrix}
	{\bs K}_m^{-1}(\hat{\bs x})
	\begin{bmatrix}
  	\partial_{x} v_m(\hat{\bs x}) \\
  	\partial_{y} v_m(\hat{\bs x})
  	\end{bmatrix}a_m(\hat{\bs x})\d\hat{\bs x}.
\end{equation}
Thus, the entries of the mass matrix and the stiffness 
matrix can be computed by means of \eqref{eq:mass} and 
\eqref{eq:stiffness} using appropriately chosen quadrature 
formulae that provide sufficient accuracy.

\subsection{Computing the right-hand side}
The right-hand side $f$ in \eqref{eq:PDE} corresponds to 
white noise. The associated discrete right-hand side ${\bs f}$ 
is thus normally distributed in accordance with ${\bs f}\sim
\mathcal{N}({\bs 0},{\bs M})$, compare \cite{Chris}. The simulation
of the desired Gaussian random field requires therefore the application 
of the matrix square root of ${\bs M}$, as a specific realization is 
generated by ${\bs f} = \sqrt{\bs M}{\bs y}$ with ${\bs y}$ being 
a sequence of i.i.d.\ $\mathcal{N}(0,1)$-distributed random variables.

The computation of the square root $\sqrt{\bs M}$ to the 
mass matrix ${\bs M}$ can be carried out in accordance with
\cite[Eq.~(4.4) and comments below]{Hale}: Let $0<m<\min
\lambda({\bs M})$ be a lower bound for smallest eigenvalue of 
the mass matrix ${\bs M}$ and $M > \max \lambda({\bs M})$ 
be an upper bound for its largest eigenvalue. Then, for some 
$\hat{K}\in\mathbb{N}$ and $\varkappa_{\bs M} := \nicefrac{M}{m}$, 
we compute
\begin{equation}\label{eq:sqrt}
  {\bs f} = \sqrt{\bs M}{\bs y}
  \approx \frac{2 E \sqrt{m}}{\pi \hat{K}}{\bs M}
  \sum_{k=1}^{\hat{K}} \frac{\dn(t_k|1-\varkappa_{\bs M}^{-1})}{\cn^2(t_k|1-\varkappa_{\bs M}^{-1})}
  	({\bs M}+w_k^2{\bs I})^{-1}{\bs y}.
\end{equation}
Here, $\sn$, $\cn$ and $\dn$ are the Jacobian elliptic 
functions (see \cite[Chpt.~16]{Abramowitz1964}), $E$ is the complete 
elliptic integral of the second kind associated with the parameter 
$\varkappa_{\bs M}^{-1}$ (see \cite[Chpt.~17]{Abramowitz1964}), and 
\[
w_k := \sqrt{m} \frac{
	\sn\left(t_k | 
	1 - \varkappa_{\bs M}^{-1} \right)}{
	\cn\left(t_k | 
	1 - \varkappa_{\bs M}^{-1} \right)},
\quad 
t_k := 
\frac{\bigl(k-\tfrac{1}{2} \bigr) E}{\hat{K}},
\quad
k\in\{1,\ldots, \hat{K} \}.
\]
Note that the computation of the approximation \eqref{eq:sqrt} 
requires only the repeated solution of (well-conditioned) linear 
systems of equations involving the mass matrix. The convergence 
towards the exact result $\sqrt{\bs M}{\bs y}$ is exponential in the 
parameter $\hat{K}$, compare \cite[Thm 4.1]{Hale}.

We remark that white noise satisfies $f\in L_{\mathbb{P}}^2(\Omega)
\otimes H^{-r}(\Gamma)$ for any $r>1$, see \cite{Kirchner,Chris} for 
example. Hence, the Gaussian random field $u$ given by \eqref{eq:PDE} 
satisfies $u\in L_{\mathbb{P}}^2(\Omega)\otimes H^{2\beta-r}(\Gamma)$. 
If $\beta>\frac{3}{2}$, linear finite elements are no more optimal and 
higher order ansatz functions pay off. More precisely, we have the 
error estimate
\[
 \|u-u_j\|_{L_{\mathbb{P}}^2(\Omega)\otimes L^2(\Gamma)}
 \lesssim 2^{-j\min\{p+1,2\beta-1-\delta\}}\quad
 \text{for any $\delta > 0$}
\]
for the mean-square error of the approximate Gaussian 
random field $u_j$ provided that the surface $\Gamma$ 
is smooth enough, compare \cite{Chris}.

\subsection{Solving equations with fractional powers}
\label{sct:fractional}
Once we have the stiffness matrix ${\bs S}$, the
mass matrix ${\bs M}$, and the right-hand side ${\bs f}$
at hand, we shall consider the discretization of the equation
\eqref{eq:PDE}. If $\beta = 1$, this is not a problem, as 
we can solve
\[
  (\kappa^2{\bs M}+{\bs S}) {\bs u} = {\bs f}
\]
in a straightforward manner. If $\beta\in\mathbb{N}$ with $\beta\ge 2$,
we successively solve equation \eqref{eq:simple solve} just 
$\beta$ times, i.e.
\[
  (\kappa^2{\bs M}+{\bs S}) {\bs w}_1 = {\bs f}\quad\text{and}\quad
   (\kappa^2{\bs M}+{\bs S}) {\bs w}_k = {\bs M}{\bs w}_{k-1}
   \ \text{for}\ k=2,\ldots,\beta,
\]
to get the approximation ${\bf u}\isdef {\bs w}_{\beta}$.
If $\beta\not\in\mathbb{N}$, we can follow 
\cite{Chris} and split $\beta$ into its integer part and the remainder 
$\in (0,1)$. We hence can assume without loss of generality 
that $\beta\in (0,1)$ since the integer part can be treated as 
described above.

The solution of \eqref{eq:PDE} for $\beta\in (0,1)$
will be performed as proposed in \cite{Bonito,Chris} by 
using the representation \eqref{eq:integral}. We choose 
$K\in\mathbb{N}$ and compute
\begin{equation}\label{eq:beta-approx}
  {\bs u}\approx {\bs u}_K = \frac{2\sin(\pi\beta)}{\sqrt{K}\pi}
  	\sum_{k=-K}^K e^{2\beta t_k}({\bs I}+e^{2 t_k}{\bs A})^{-1} {\bs f}
\end{equation}
with quadrature points $t_k = k/\sqrt{K}$ and ${\bs A} = 
\kappa^2{\bs M}+{\bs S}$. As shown in \cite[Lemma 2.3]{Chris}, 
we obtain exponential convergence with respect to $\sqrt{K}$ 
according to
\begin{equation}\label{eq:K-convergence}
  \|{\bs u}-{\bs u}_K\|\lesssim e^{-2\min\{\beta,1-\beta\}\sqrt{K}}\|{\bs f}\|,
\end{equation}
where the constant involved does not depend on $\beta$.

The convergence rate \eqref{eq:K-convergence} of the sinc 
quadrature depends crucially on $\beta$. It is best for $\beta =
\nicefrac{1}{2}$ and gets significantly worse when $\beta$ is 
close to $0$ or $1$. Thus, we propose the following modification 
compared to \cite{Chris} if $\beta$ is close to an integer $n$:
\begin{itemize}
\item
If $\beta = n+\beta'$ with $\beta'\in (0,\nicefrac{1}{3})$, we split 
$\beta = n-1 + 2\beta''$ with $\beta''=\nicefrac{(\beta'+1)}{2}\in 
(\nicefrac{1}{2},\nicefrac{2}{3})$ instead of using the integer and 
fractional parts of $\beta$.
\item
If $\beta = n+\beta'$ with $\beta'\in (\nicefrac{2}{3},1)$, we split 
$\beta = n + 2\beta''$ with $\beta''=\nicefrac{\beta'}{2}\in (\nicefrac{1}{3},
\nicefrac{1}{2})$ instead of using the integer and fractional parts
 of $\beta$. 
\end{itemize}
One readily veryfies from \eqref{eq:K-convergence} that 
the convergence for $\beta''$ is in both cases faster than 
for $\beta'$. Although we have to apply \eqref{eq:beta-approx}
twice in this new approach, we obtain improved computation 
times in comparison to the original approach of \cite{Chris} 
as we demonstrate in Section~\ref{sec:numerix}. 


\subsection{Preconditioning}\label{sec:solver}
In isogeometric analysis, we naturally obtain a sequence of 
nested ansatz spaces as the step size $h$ tends to zero. This 
allows for geometric multigrid techniques for the solution of the 
linear system of equations resulting from the discretization of 
the variational formulation \eqref{eq:var_form}. In our numerical 
experiments, we consider both, the multiplicative multigrid method 
by means of the V-cycle and the additive multigrid method by means 
of the BPX preconditioner. These methods have been introduced in 
the context isogeometric analysis in \cite{BPX,MGM}. Thus, we 
obtain a linear cost-complexity rate for the computation of the 
numerical approximation to \eqref{eq:var_form}. Nonetheless, 
a polynomial dependency of the complexity on the polynomial 
degree $p$ remains.

\section{Numerical experiments}
\label{sec:numerix}
\subsection{Results for the sphere}
\begin{figure}[hbt]
  \centering
  \includegraphics[width=0.6\textwidth]{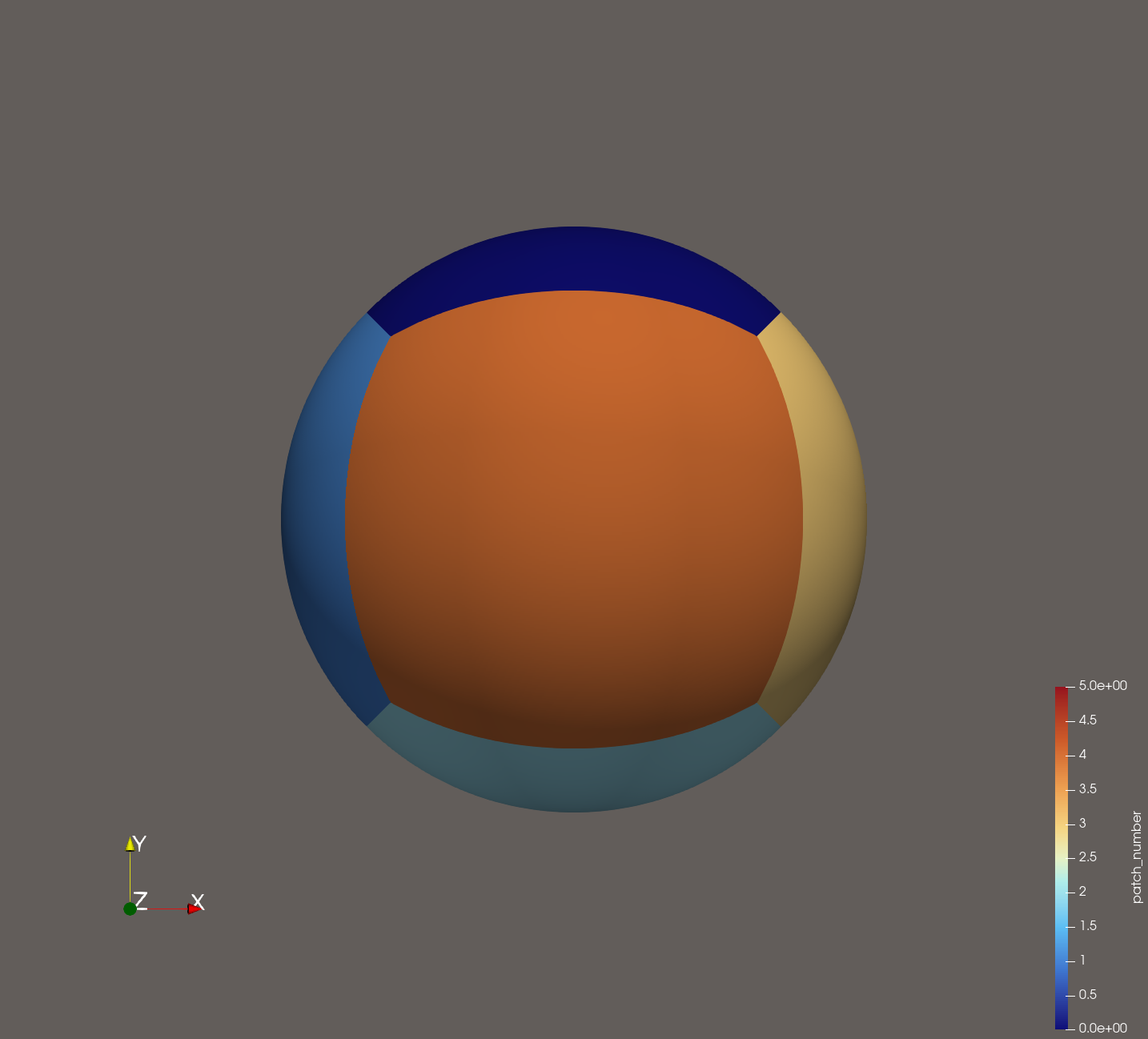}
  \caption{\label{fig:sphere}The unit sphere, consisting of 6 patches.}
\end{figure}

The first surface we consider is the unit sphere $\mathbb{S}^2 = 
\{(x,y,z)\in \mathbb{R}^3: x^2 + y^2 + z^2 = 1\}$, which is parametrized 
by six patches as seen in Figure~\ref{fig:sphere}. In order to test our 
implementation, we use the spherical harmonics $Y_{\ell,m}$, as they 
are known to be the eigenfunctions of the Laplace-Beltrami operator 
on the unit sphere, see \cite{atkinson_spherical_2012} for example. 
In our experiments, we especially study the influence of the 
polynomial degree $p$ and the refinement level $j$ on the
proposed algorithms.

\subsubsection{Rate of convergence}
\begin{figure}[hbt]
\begin{center}
\begin{tikzpicture}[scale=1]
\begin{semilogyaxis}[
xlabel={refinement level $j$},
ylabel={$L^2$-error},
grid=major,
xtick={0,1,2,3,4,5,6,7,8,9},
ytick={1,1e-2,1e-4,1e-6,1e-8,1e-10,1e-12},
xmin=0,
xmax=9,
ymin=1e-13,
ymax=1,
legend style={
at={(1.05,0.5)}, 
anchor=west 
}
]
\addlegendentry{$p=1$}
\addplot[blue, mark=triangle*, mark options={scale=1}] table[x index=0, y index=1]{
0  0.120221
1  0.0303486
2  0.00763282
3  0.00190615
4  0.000476429
5  0.000119108
6  2.97778E-05
7  7.44456E-06
8  1.86115E-06
9  4.65288E-07
};
\addlegendentry{$p=2$}
\addplot[green, mark=*, mark options={scale=1}] table[x index=0, y index=1]{
0  0.0234786
1  0.00495968
2  0.000525151
3  5.76864E-05
4  6.97993E-06
5  8.65378E-07
6  1.0795E-07
7  1.34868E-08
8  1.68563E-09
9  2.10697E-10
};
\addlegendentry{$p=3$}
\addplot[red, mark=pentagon*, mark options={scale=1}] table[x index=0, y index=1]{
0  0.00441012
1  0.000865953
2  8.10752E-05
3  3.73186E-06
4  2.16352E-07
5  1.33029E-08
6  8.29186E-10
7  5.18254E-11
8  3.2433E-12
9  2.24224E-13
};
\addlegendentry{$p=4$}
\addplot[lila, mark=diamond*, mark options={scale=1}] table[x index=0, y index=1]{
0  0.000598262
1  0.000159845
2  1.73803E-05
3  3.32503E-07
4  8.75729E-09
5  2.61405E-10
6  8.07666E-12
7  2.59847E-13
};
\addlegendentry{$p=5$}
\addplot[orange!80!black, mark=square*, mark options={scale=1}] table[x index=0, y index=1]{
0  0.000217381
1  2.82315E-05
2  3.86E-06
3  3.7252E-08
4  4.35305E-10
5  6.26225E-12
6  1.77647E-13
};
\addlegendentry{$2^{-2j}$}
\addplot[blue, dashed] table[x index=0, y index=1]
{
0   0.25
9   9.5367e-07
};
\addlegendentry{$2^{-3j}$}
\addplot[green, dashed] table[x index=0, y index=1]
{
0   0.05
9   3.7253e-10
};
\addlegendentry{$2^{-4j}$}
\addplot[red, dashed] table[x index=0, y index=1]
{
0   0.025
9   3.6380e-13
};
\addlegendentry{$2^{-5j}$}
\addplot[lila, dashed] table[x index=0, y index=1]
{
0   0.025
9   7.1054e-16
};
\addlegendentry{$2^{-6j}$}
\addplot[orange!80!black, dashed] table[x index=0, y index=1]
{
0   0.025
9   1.3878e-18
};
\end{semilogyaxis}
\end{tikzpicture}
\end{center}
  \caption{\label{fig:convergence}Convergence rate of the approximate solution
  to \eqref{eq:sphere} in dependence of the refinement level $j$ and the 
  polynomial degree $p$.}
\end{figure}
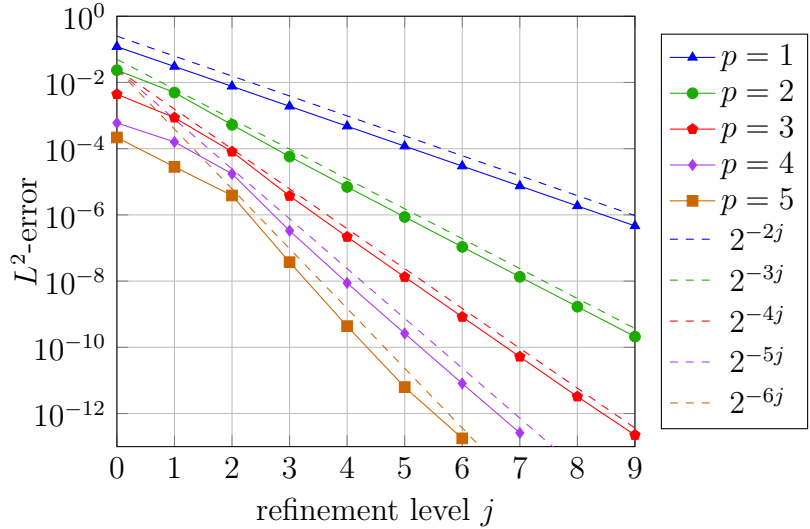

The isogeometric finite element implementation of the 
Laplace-Beltrami operator is tested on the unit sphere 
with the problem 
\begin{equation}\label{eq:sphere}
  (\kappa^2-\Delta_\Gamma) u = Y_{1,-1}
  \quad\text{on $\Gamma = \mathbb{S}^2$},
\end{equation}
which has the analytical solution 
\[
  u = \frac{1}{2+\kappa^2} Y_{1,-1}.
\]
We choose $\kappa=1$ as a convenient value, but the 
choice of $\kappa$ does not impact the results. We discretize 
the equation \eqref{eq:sphere} as described in the previous 
sections and solve the resulting linear system of equations 
by using Eigen's conjugate gradient method with tolerance 
$\varepsilon = 10^{-12}$ for different refinement levels $j$ 
and corresponding mesh sizes $h_j \sim 2^{-j}$. The errors 
obtained for different polynomial degrees $1\leq p\leq 5$ of 
our ansatz space are reported in Figure~\ref{fig:convergence}.
As one figures out thereof, the convergence rates observed 
match the theoretical rates from \eqref{eq:approximation}:
\[
  \|u-u_j\|_{L^2(\mathbb{S}^2)} 
  	\lesssim 2^{-j(p+1)} \|u\|_{H^{p+1}(\mathbb{S}^2)}.
\]

\subsubsection{Preconditioning}
\begin{figure}[hbt]
\begin{center}
\begin{tikzpicture}[scale=0.7]
\begin{semilogyaxis}[
title = {BPX with diagonal scaling},
xlabel={Refinement level $j$},
ylabel={Time in seconds},
xtick={0,1,2,3,4,5,6,7,8,9,10},
ytick={1,1e-2,1e-4,1e-6,1e-8,1e-10,1e-12},
xmin=0,
xmax=9,
ymin=1e-7,
ymax=500,
legend style={
at={(0.7,0.02)}, 
anchor=south west 
}
]
\addlegendentry{$p = 1$}
\addplot[blue, mark=triangle*, mark options={scale=1}] table[x index=0, y index=1]{
0 2.067E-06
1 6.889E-06
2 4.1659E-05
3 0.00022773
4 0.00105291
5 0.00454129
6 0.0208683
7 0.0870297
8 0.562899
9 2.90975
};
\addlegendentry{$p = 2$}
\addplot[green, mark=*, mark options={scale=1}] table[x index=0, y index=1]{
0 6.333E-06
1 2.2653E-05
2 0.000136589
3 0.00102633
4 0.00482397
5 0.0180616
6 0.0764475
7 0.591945
8 2.66914
9 10.4082
};
\addlegendentry{$p = 3$}
\addplot[red, mark=pentagon*, mark options={scale=1}] table[x index=0, y index=1]{
0 2.8186E-05
1 9.8515E-05
2 0.000581392
3 0.00424015
4 0.0185709
5 0.0711165
6 0.282141
7 1.80796
8 8.19857
9 34.5058
};
\addlegendentry{$p = 4$}
\addplot[lila, mark=diamond*, mark options={scale=1}] table[x index=0, y index=1]{
0 4.5713E-05
1 0.000324766
2 0.0017793
3 0.0120704
4 0.0565372
5 0.225701
6 0.932251
7 5.64073
8 24.134
9 107.047
};
\addlegendentry{$p = 5$}
\addplot[orange!80!black, mark=square*, mark options={scale=1}] table[x index=0, y index=1]{
0 0.000126229
1 0.00145773
2 0.00671603
3 0.0378621
4 0.18756
5 0.69601
6 3.37794
7 16.7946
8 71.33
9 302.312
};
\addlegendentry{$4^j$}
\addplot[dashed] table[x index=0, y index=1]
{
 0  2.5000e-06
 1  1.0000e-05
 2  4.0000e-05
 3  1.6000e-04
 4  6.4000e-04
 5  2.5600e-03
 6  1.0240e-02
 7  4.0960e-02
 8  1.6384e-01
 9  6.5536e-01
};
\end{semilogyaxis}
\end{tikzpicture}
\begin{tikzpicture}[scale=0.7]
\begin{semilogyaxis}[
title = {BPX with SSOR},
xlabel={Refinement level $j$},
ylabel={Time in seconds},
xtick=\empty, 
xtick={0,1,2,3,4,5,6,7,8,9,10},
ytick={1,1e-2,1e-4,1e-6,1e-8,1e-10,1e-12},
xmin=0,
xmax=9,
ymin=1e-7,
ymax=500,
legend style={
at={(0.7,0.02)}, 
anchor=south west 
}
]
\addlegendentry{$p = 1$}
\addplot[blue, mark=triangle*, mark options={scale=1}] table[x index=0, y index=1]{
0 5.679E-06
1 2.9634E-05
2 0.000120754
3 0.000487688
4 0.00225914
5 0.0102894
6 0.0410177
7 0.168223
8 0.946475
9 4.68366
};
\addlegendentry{$p = 2$}
\addplot[green, mark=*, mark options={scale=1}] table[x index=0, y index=1]{
0 1.755E-05
1 0.000104914
2 0.000414928
3 0.00149175
4 0.00577575
5 0.0225448
6 0.0912088
7 0.456226
8 2.41012
9 10.6912
};
\addlegendentry{$p = 3$}
\addplot[red, mark=pentagon*, mark options={scale=1}] table[x index=0, y index=1]{
0 7.9063E-05
1 0.000601131
2 0.00189032
3 0.00599839
4 0.0194262
5 0.0618465
6 0.248319
7 1.41202
8 6.48655
9 28.6615
};
\addlegendentry{$p = 4$}
\addplot[lila, mark=diamond*, mark options={scale=1}] table[x index=0, y index=1]{
0 0.00046348
1 0.00216485
2 0.00662735
3 0.0179772
4 0.05475
5 0.186781
6 0.727077
7 4.13266
8 17.9312
9 78.2186
};
\addlegendentry{$p = 5$}
\addplot[orange!80!black, mark=square*, mark options={scale=1}] table[x index=0, y index=1]{
0 0.00167788
1 0.00791104
2 0.0207667
3 0.052159
4 0.156432
5 0.511468
6 2.4148
7 11.5107
8 48.7675
9 203.325
};
\addlegendentry{$4^j$}
\addplot[dashed] table[x index=0, y index=1]
{
0   5.0000e-06
1   2.0000e-05
2   8.0000e-05
3   3.2000e-04
4   1.2800e-03
5   5.1200e-03
6   2.0480e-02
7   8.1920e-02
8   3.2768e-01 
9   1.3107e+00
};
\end{semilogyaxis}
\end{tikzpicture}
\caption{\label{fig:sphere_BPX}Runtime of the conjugate gradient 
with BPX-preconditioning in case of a diagonal scaling on each level 
(left) and in case of an SSOR preconditioning on each level (right).}
\end{center}
\end{figure}
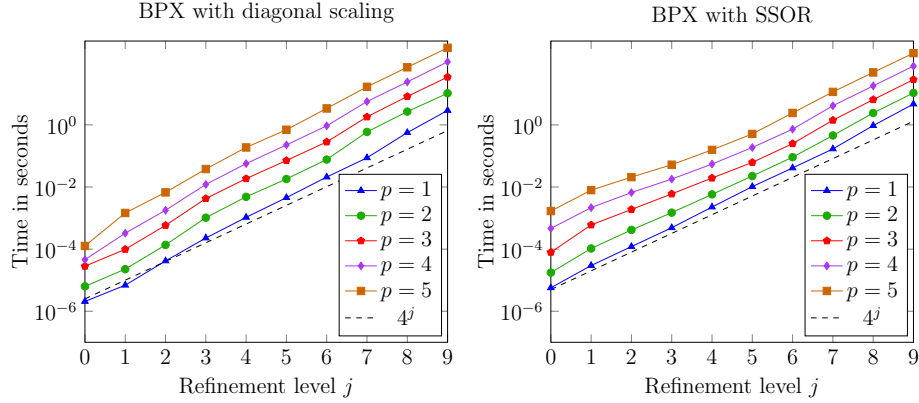

We shall next study the preconditioning of the linear system of 
equations which arises from the isogeometric finite element 
discretization of the boundary value problem \eqref{eq:sphere}. 
To this end, we apply the conjugate gradient method, which we 
precondition by the BPX-preconditioner with either a diagonal 
scaling on each level (``BPX-diag'') or an SSOR preconditioning 
on each level (``BPX-SSOR''), cf.~\cite{BPX} for the details. The
computational complexity is then linear in the number of degrees 
of freedom if the system of linear equations are solved with 
fixed accuracy, but it depends on the polynomial degree $p$. 
This has already been reported in \cite{BPX} and matches the 
theory since the spectral condition number of the preconditioned 
system matrix is independent of the mesh size $h_j$ but 
depends on $p$. Indeed, we observe a basically linear increase 
of the runtimes in Figure~\ref{fig:sphere_BPX} in the number of 
degrees of freedom for fixed polynomial degree $p$. As it can be 
expected, the BPX-SSOR variant takes less iterations per refinement 
level $j$, compare Table\ \ref{tab:sphere_BPX}. However, it is also 
computationally more expensive. Hence, we barely notice a 
difference in the runtimes in practice between both variants 
of the BPX-preconditioner.

\begin{table}[hbt]
\begin{center}
{\small
\begin{tabular}{|c|ccccc|ccccc|}\cline{2-11}
\multicolumn{1}{c|}{}
& \multicolumn{5}{c|}{BPX with diagonal scaling} & \multicolumn{5}{c|}{BPX with SSOR}\\\hline
$j$ & $p=1$ & $p=2$ & $p=3$ & $p=4$ & $p=5$ &$p=1$ & $p=2$ & $p=3$ & $p=4$ & $p=5$ \\\hline
0 & 1 & 2 & 4 & 9 & 15 & 4 & 8 & 19 & 36 & 70 \\
1 & 2 & 4 & 9 & 15 & 39 & 11 & 16 & 29 & 52 & 105 \\
2 & 8 & 11 & 24 & 34 & 76 & 18 & 21 & 34 & 64 & 115 \\
3 & 17 & 28 & 62 & 93 & 177 & 21 & 25 & 38 & 68 & 122 \\
4 & 21 & 39 & 87 & 163 & 343 & 23 & 27 & 42 & 73 & 135 \\
5 & 24 & 44 & 99 & 204 & 423 & 26 & 29 & 44 & 77 & 141 \\
6 & 27 & 49 & 108 & 223 & 462 & 28& 31 & 47 & 81 & 147 \\
7 & 29 & 53 & 116 & 238 & 494 & 29 & 33 & 49 & 85 & 152 \\
8 & 30 & 56 & 124 & 253 & 525 & 31 & 34 & 51 & 88 & 157 \\
9 & 32 & 59 & 130 & 264 & 548 & 32 & 35 & 53 & 91 & 161 \\\hline
\end{tabular}}
\caption{\label{tab:sphere_BPX}Number of iterations of the preconditioned 
conjugate gradient with BPX-preconditioning in case of a diagonal scaling 
on each level and in case of an SSOR preconditioning on each level.}
\end{center}
\end{table}

We like to mention that we also tested a multiplicative multigrid
using the $V$-cycle. It turned out to be less robust with respect 
to higher polynomial degree $p$ and with respect to the geometry.
In all the test cases we considered, the conjugate gradient method 
in combination with the BPX-preconditioner was more efficient.

\subsubsection{Sinc quadrature}
We next consider the fractional problem
\[
  (\kappa^2-\Delta_\Gamma)^\beta u = Y_{1,-1}
  \quad\text{on $\Gamma = \mathbb{S}^2$}
\]
with analytical solution
\[
u = \left(\frac{1}{2+\kappa^2}\right)^{\beta} Y_{1,-1}.
\]
We solve this fractional problem by using the sinc quadrature for 
$p=3$ and $j=5$ (where we can expect an discretization error of
about $10^{-8})$ and the number of quadrature points $K$ ranging 
from $2$ to $500$. The convergence rate we observe matches the 
predicted rate of $\operatorname{exp}(C_\beta \sqrt{K})$, where 
$C_\beta = -2\min\{\beta, 1-\beta\}$, compare \eqref{eq:K-convergence}. 
It should also be noted that the multiplicative constants hidden 
in the error rates seem to be independent of $\beta$, as is 
important for our new proposed method of solving \eqref{eq:PDE} 
for fractional $\beta$ being close to an integer.

\begin{figure}[hbt]
\begin{center}
\begin{tikzpicture}[scale=1]
\begin{semilogyaxis}[
xlabel={Number $K$ of quadrature points},
ylabel={$L^2$-error},
xtick=\empty, 
extra x ticks={1.4142,2.2361,3.1623,4.4721,7.0711,10,14.1421,22.3607},
extra x tick labels={$2\ \ $,$5\;$,$\;10$,$\quad20$,50,100,200,500},
ytick={1,1e-2,1e-4,1e-6,1e-8,1e-10,1e-12},
xmin=0,
xmax=25,
ymin=1e-8,
ymax=2,
legend style={
at={(1.05,0.5)}, 
anchor=west 
}
]
\addlegendentry{$\beta = 0.15$}
\addplot[blue, mark=triangle*, mark options={scale=1}] table[x index=1, y index=2]{
2  1.4142  1.1257
5  2.2361  0.919747
10  3.1623  0.711373
20  4.4721  0.487031
50  7.0711  0.226124
100  10.00  0.0945236
200  14.1421  0.0274567
500  22.3607  0.00288245
};
\addlegendentry{$\beta = 0.30$}
\addplot[green, mark=*, mark options={scale=1}] table[x index=1, y index=2]{
2  1.4142  0.599642
5  2.2361  0.39712
10  3.1623  0.235926
20  4.4721  0.110019
50  7.0711  0.0236488
100  10.00  0.00412929
200  14.1421  0.000347025
500  22.3607  2.52484E-06
};
\addlegendentry{$\beta = 0.50$}
\addplot[red, mark=pentagon*, mark options={scale=1}] table[x index=1, y index=2]{
2  1.4142  0.278578
5  2.2361  0.143127
10  3.1623  0.0610361
20  4.4721  0.0173034
50  7.0711  0.00134231
100  10.00  7.32839E-05
200  14.1421  1.18205E-06
500  22.3607  2.30435E-08
};
\addlegendentry{$\beta = 0.70$}
\addplot[lila, mark=diamond*, mark options={scale=1}] table[x index=1, y index=2]{
2  1.4142  0.254041
5  2.2361  0.153293
10  3.1623  0.0849038
20  4.4721  0.0377388
50  7.0711  0.00791263
100  10.00  0.00137694
200  14.1421  0.000115677
500  22.3607  8.41762E-07
};
\addlegendentry{$\beta = 0.85$}
\addplot[orange!80!black, mark=square*, mark options={scale=1}] table[x index=1, y index=2]{
2  1.4142  0.392153
5  2.2361  0.311779
10  3.1623  0.238279
20  4.4721  0.162469
50  7.0711  0.0753706
100  10.00  0.0314993
200  14.1421  0.00913155
500  22.3607  0.000778864
};
\addlegendentry{$\exp(C_{0.15}\sqrt{K})$}
\addplot[blue, dashed] table[x index=1, y index=2]
{
2  1.4142  0.75
500  22.3607  1.6e-3
};
\addlegendentry{$\exp(C_{0.30}\sqrt{K})$}
\addplot[green, dashed] table[x index=1, y index=2]
{
2  1.4142  0.4
500  22.3607  1.5e-6
};
\addlegendentry{$\exp(C_{0.50}\sqrt{K})$}
\addplot[red, dashed] table[x index=1, y index=2]
{
2  1.4142  0.5
500  22.3607  5.9488e-10
};
\end{semilogyaxis}
\end{tikzpicture}
\caption{Convergence of the sinc quadrature with respect 
to the number $K$ of quadrature points for fixed level of 
refinement $j=6$ and polynomial degree $p = 3$.}
\end{center}
\end{figure}
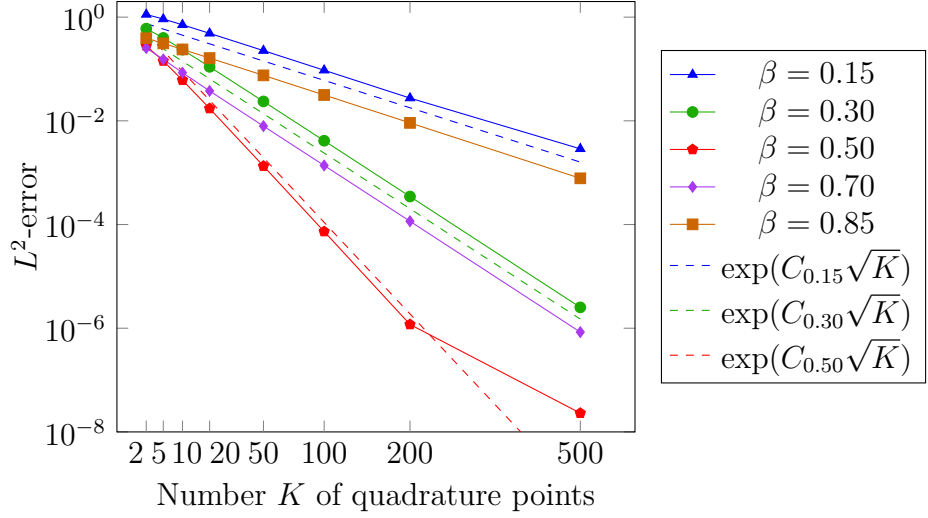

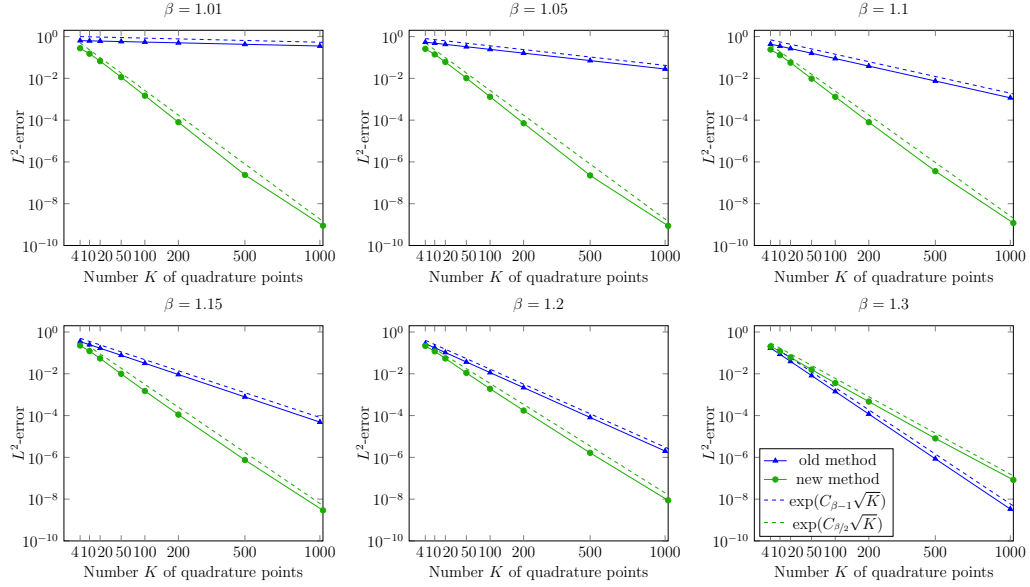
\begin{figure}[hbt]
\begin{center}
\begin{tikzpicture}[scale=0.5]
\begin{semilogyaxis}[
title = {$\beta = 1.01$},
xlabel={Number $K$ of quadrature points},
ylabel={$L^2$-error},
xtick=\empty, 
extra x ticks={2,3.1623,4.4721,7.0711,10,14.1421,22.3607,31.6228},
extra x tick labels={$4\ \ $,$10\,$,$\ \ 20$,$\ 50$,100,200,500,1000},
ytick={1,1e-2,1e-4,1e-6,1e-8,1e-10,1e-12},
xmin=0,
xmax=32,
ymin=1e-10,
ymax=2,
]
\addplot[blue, mark=triangle*, mark options={scale=1}] table[x index=1, y index=2]{
4  2  0.639882
10  3.1623  0.625621
20  4.4721  0.609506
50  7.0711  0.578649
100  10.00  0.545729
200  14.1421  0.502342
500  22.3607  0.426201
1000  31.6228  0.354132
};
\addplot[green, mark=*, mark options={scale=1}] table[x index=1, y index=2]{
4  2  0.277807
10  3.1623  0.1517
20  4.4721  0.0668502
50  7.0711  0.0114743
100  10.00  0.00148281
200  14.1421  8E-05
500  22.3607  2.39E-07
1000  32  8.87E-10
};
\addplot[blue, dashed] table[x index=1, y index=2]
{
4  2  1
1225  35  0.5
};
\addplot[green, dashed] table[x index=1, y index=2]
{
4  2  0.5
1225  35  0.2E-9
};
\end{semilogyaxis}
\end{tikzpicture}
\begin{tikzpicture}[scale=0.5]
\begin{semilogyaxis}[
title = {$\beta = 1.05$},
xlabel={Number $K$ of quadrature points},
ylabel={$L^2$-error},
xtick=\empty, 
extra x ticks={2,3.1623,4.4721,7.0711,10,14.1421,22.3607,31.6228},
extra x tick labels={$4\ \ $,$10\,$,$\ \ 20$,$\ 50$,100,200,500,1000},
ytick={1,1e-2,1e-4,1e-6,1e-8,1e-10,1e-12},
xmin=0,
xmax=32,
ymin=1e-10,
ymax=2,
]
\addplot[blue, mark=triangle*, mark options={scale=1}] table[x index=1, y index=2]{
4  2  0.541363
10  3.1623  0.483606
20  4.4721  0.424444
50  7.0711  0.327348
100  10.00  0.244242
200  14.1421  0.161411
500  22.3607  0.0709592
1000  31.6228  0.0281036
};
\addplot[green, mark=*, mark options={scale=1}] table[x index=1, y index=2]{
4  2  0.259498
10  3.1623  0.140157
20  4.4721  0.0610391
50  7.0711  0.0102924
100  10.00  0.00131604
200  14.1421  7.12E-05
500  22.3607  2.26E-07
1000  32  8.7E-10
};
\addplot[blue, dashed] table[x index=1, y index=2]
{
4  2  0.8
1225  35  0.03
};
\addplot[green, dashed] table[x index=1, y index=2]
{
4  2  0.5
1225  35  0.2E-9
};
\end{semilogyaxis}
\end{tikzpicture}
\begin{tikzpicture}[scale=0.5]
\begin{semilogyaxis}[
title = {$\beta = 1.1$},
xlabel={Number $K$ of quadrature points},
ylabel={$L^2$-error},
xtick=\empty, 
extra x ticks={2,3.1623,4.4721,7.0711,10,14.1421,22.3607,31.6228},
extra x tick labels={$4\ \ $,$10\,$,$\ \ 20$,$\ 50$,100,200,500,1000},
ytick={1,1e-2,1e-4,1e-6,1e-8,1e-10,1e-12},
xmin=0,
xmax=32,
ymin=1e-10,
ymax=2,
]
\addplot[blue, mark=triangle*, mark options={scale=1}] table[x index=1, y index=2]{
4  2  0.436195
10  3.1623  0.347938
20  4.4721  0.267996
50  7.0711  0.159408
100  10.00  0.088743
200  14.1421  0.0387582
500  22.3607  0.00749057
1000  31.6228  0.00117495
};
\addplot[green, mark=*, mark options={scale=1}] table[x index=1, y index=2]{
4  2  0.240756
10  3.1623  0.12941
20  4.4721  0.056349
50  7.0711  0.00968722
100  10.00  0.00130519
200  14.1421  7.98E-05
500  22.3607  3.61E-07
1000  32  1.19E-09
};
\addplot[blue, dashed] table[x index=1, y index=2]
{
4  2  0.7
1225  35  0.001
};
\addplot[green, dashed] table[x index=1, y index=2]
{
4  2  0.4
1225  35  0.3E-9
};
\end{semilogyaxis}
\end{tikzpicture}
\begin{tikzpicture}[scale=0.5]
\begin{semilogyaxis}[
title = {$\beta = 1.15$},
xlabel={Number $K$ of quadrature points},
ylabel={$L^2$-error},
xtick=\empty, 
extra x ticks={2,3.1623,4.4721,7.0711,10,14.1421,22.3607,31.6228},
extra x tick labels={$4\ \ $,$10\,$,$\ \ 20$,$\ 50$,100,200,500,1000},
ytick={1,1e-2,1e-4,1e-6,1e-8,1e-10,1e-12},
xmin=0,
xmax=32,
ymin=1e-10,
ymax=2,
]
\addplot[blue, mark=triangle*, mark options={scale=1}] table[x index=1, y index=2]{
4  2  0.348839
10  3.1623  0.24829
20  4.4721  0.167802
50  7.0711  0.0769764
100  10.00  0.0319741
200  14.1421  0.00922879
500  22.3607  0.000784104
1000  31.6228  4.87E-05
};
\addplot[green, mark=*, mark options={scale=1}] table[x index=1, y index=2]{
4  2  0.226242
10  3.1623  0.122353
20  4.4721  0.0542094
50  7.0711  0.00995245
100  10.00  0.0015009
200  14.1421  0.000111583
500  22.3607  7.46E-07
1000  32  2.91E-09
};
\addplot[blue, dashed] table[x index=1, y index=2]
{
4  2  0.5
1225  35  3E-05
};
\addplot[green, dashed] table[x index=1, y index=2]
{
4  2  0.4
1225  35  0.8E-9
};
\end{semilogyaxis}
\end{tikzpicture}
\begin{tikzpicture}[scale=0.5]
\begin{semilogyaxis}[
title = {$\beta = 1.2$},
xlabel={Number $K$ of quadrature points},
ylabel={$L^2$-error},
xtick=\empty, 
extra x ticks={2,3.1623,4.4721,7.0711,10,14.1421,22.3607,31.6228},
extra x tick labels={$4\ \ $,$10\,$,$\ \ 20$,$\ 50$,100,200,500,1000},
ytick={1,1e-2,1e-4,1e-6,1e-8,1e-10,1e-12},
xmin=0,
xmax=32,
ymin=1e-10,
ymax=2
]
\addplot[blue, mark=triangle*, mark options={scale=1}] table[x index=1, y index=2]{
4  2  0.277068
10  3.1623  0.17577
20  4.4721  0.104181
50  7.0711  0.0368521
100  10.00  0.0114214
200  14.1421  0.00217864
500  22.3607  8.14E-05
1000  31.6228  2E-06
};
\addplot[green, mark=*, mark options={scale=1}] table[x index=1, y index=2]{
4  2  0.215556
10  3.1623  0.118603
20  4.4721  0.0543666
50  7.0711  0.011048
100  10.00  0.00191843
200  14.1421  0.00017422
500  22.3607  1.64E-06
1000  32  8.76E-09
};
\addplot[blue, dashed] table[x index=1, y index=2]
{
4  2  0.4
1225  35  0.8E-6
};
\addplot[green, dashed] table[x index=1, y index=2]
{
4  2  0.3
1225  35  3E-9
};
\end{semilogyaxis}
\end{tikzpicture}
\begin{tikzpicture}[scale=0.5]
\begin{semilogyaxis}[
title = {$\beta = 1.3$},
xlabel={Number $K$ of quadrature points},
ylabel={$L^2$-error},
xtick=\empty, 
extra x ticks={2,3.1623,4.4721,7.0711,10,14.1421,22.3607,31.6228},
extra x tick labels={$4\ \ $,$10\,$,$\ \ 20$,$\ 50$,100,200,500,1000},
ytick={1,1e-2,1e-4,1e-6,1e-8,1e-10,1e-12},
xmin=0,
xmax=32,
ymin=1e-10,
ymax=2,
legend style={
at={(0.02,0.02)}, 
anchor=south west
}
]
\addlegendentry{old method}
\addplot[blue, mark=triangle*, mark options={scale=1}] table[x index=1, y index=2]{
4  2  0.172231
10  3.1623  0.0863101
20  4.4721  0.0391939
50  7.0711  0.00822298
100  10.00  0.00141824
200  14.1421  0.000118152
500  22.3607  8.53E-07
1000  31.6228  3.35E-09
};
\addlegendentry{new method}
\addplot[green, mark=*, mark options={scale=1}] table[x index=1, y index=2]{
4  2  0.204088
10  3.1623  0.1197
20  4.4721  0.0608975
50  7.0711  0.0158742
100  10.00  0.00363808
200  14.1421  0.000467412
500  22.3607  8.07E-06
1000  32  8.3E-08
};
\addlegendentry{$\exp(C_{\beta-1}\sqrt{K})$}
\addplot[blue, dashed] table[x index=1, y index=2]
{
4  2  0.25
1225  35  0.8E-9
};
\addlegendentry{$\exp(C_{\nicefrac{\beta}{2}}\sqrt{K})$}
\addplot[green, dashed] table[x index=1, y index=2]
{
4  2  0.3
1225  35  3E-8
};
\end{semilogyaxis}
\end{tikzpicture}
\caption{\label{fig:beta}Comparison of both proposed methods 
for solving the fractional operator for certain values of $\beta
\in(1,1.3]$. When $\beta\in (1,\nicefrac{5}{4})$, then solving 
twice for $\nicefrac{\beta}{2}$ is superior (green graphs) to 
solving first for $\beta' = 1$ and then for $\beta' = \beta - 1$ 
(blue graphs).}
\end{center}
\end{figure}

We shall next compare in Figure~\ref{fig:beta} the two 
methods for approximating $(\kappa^2+\Delta_\Gamma)^{-\beta} f$ 
for $1 < \beta \le\nicefrac{1}{3}$. We expect that the old method, 
where we split $\beta$ into $1$ and $\beta-1$ should be faster 
for all $\beta\in (\nicefrac{4}{3},\nicefrac{3}{2})$, whereas the 
new method, where we split $\beta$ into $\nicefrac{\beta}{2}$ 
twice, is faster for $\beta \in (1,\nicefrac{4}{3})$. We apply 
both methods with the same total number of quadrature 
nodes used, taking into account that we have to apply the 
quadrature twice for the new method. Indeed, as seen in 
Figure~\ref{fig:beta}, the new method is faster for $\beta
\le\nicefrac{5}{4}$ which is a slightly smaller range than
predicted.

\subsubsection{Matrix square root}
We test the convergence of the expansion \eqref{eq:sqrt} for 
approximating the multiplication of a given vector ${\bs y}$ 
with the square root of the mass matrix ${\bs M}$ by applying 
it twice to a random vector ${\bs y}$ and then calculating 
the error 
\[
\Big\|{\bs M}{\bs y} - {\bs M}_{\hat{K}}^{\sqrt{}}{\bs M}_{\hat{K}}^{\sqrt{}} {\bs y}\Big\|.
\] 
This gives an error which is of the same order of magnitude as 
\[
\Big\|\sqrt{\bs M} {\bs y} - {\bs M}_{\hat{K}}^{\sqrt{}} {\bs y}\Big\|
\] 
without needing to calculate the exact matrix square 
root. Note that the quantity $\varkappa_{\bs M}$ appearing
in \eqref{eq:sqrt} has been estimated by computing the largest 
and smallest eigenvalue of ${\bs M}$ by 20 steps of the 
(inverse) power method.

\begin{figure}[hbt]
\begin{center}
\begin{tikzpicture}[scale=0.7]
\begin{semilogyaxis}[
xlabel={Number $\hat{K}$ of quadrature points},
ylabel={$L^2$-error},
xtick=\empty, 
xtick={2,4,6,8,10,12,14,16,18,20},
ytick={1,1e-2,1e-4,1e-6,1e-8,1e-10,1e-12},
xmin=0,
xmax=22,
ymin=1e-16,
ymax=2,
legend style={
at={(0.01,0.01)}, 
anchor=south west 
}
]
\addlegendentry{$p=1$}
\addplot[blue, mark=triangle*, mark options={scale=1}] table[x index=0, y index=1]{
2  0.00707126
4  8.13E-06
6  9.76E-09
8  1.17E-11
10  1.4E-14
};
\addlegendentry{$p=2$}
\addplot[green, mark=*, mark options={scale=1}] table[x index=0, y index=1]{
2  0.045418
4  0.000354993
6  2.84E-06
8  2.25E-08
10  1.78E-10
12  1.42E-12
14  1.13E-14
};
\addlegendentry{$p=3$}
\addplot[red, mark=pentagon*, mark options={scale=1}] table[x index=0, y index=1]{
2  0.105092
4  0.00253605
6  5.19E-05
8  1.08E-06
10  2.28E-08
12  4.74E-10
14  9.98E-12
16  2.08E-13
18  4.58E-15
};
\addlegendentry{$p=4$}
\addplot[lila, mark=diamond*, mark options={scale=1}] table[x index=0, y index=1]{
2  0.174749
4  0.00894825
6  0.000333611
8  1.31E-05
10  5.03E-07
12  1.98E-08
14  7.61E-10
16  2.97E-11
18  1.16E-12
20  4.55E-14
};
\addlegendentry{$p=5$}
\addplot[orange!80!black, mark=square*, mark options={scale=1}] table[x index=0, y index=1]{
2  0.244923
4  0.0215028
6  0.00124307
8  7.3E-05
10  4.51E-06
12  2.61E-07
14  1.6E-08
16  9.6E-10
18  5.81E-11
20  3.45E-12
};
\addlegendentry{$p=6$}
\addplot[magenta, mark=o, mark options={scale=1}] table[x index=0, y index=1]{
2  0.340275
4  0.0396697
6  0.00350875
8  0.000282571
10  2.4E-05
12  2.05E-06
14  1.71E-07
16  1.44E-08
18  1.21E-09
20  1.03E-10
};
\end{semilogyaxis}
\end{tikzpicture}
\begin{tikzpicture}[scale=0.7]
\begin{semilogyaxis}[
xlabel={Number $\hat{K}$ of quadrature points},
ylabel={$L^2$-error},
xtick=\empty, 
xtick={2,4,6,8,10,12,14,16,18,20},
ytick={1,1e-2,1e-4,1e-6,1e-8,1e-10,1e-12},
xmin=0,
xmax=22,
ymin=1e-16,
ymax=2,
legend style={
at={(0.01,0.01)}, 
anchor=south west 
}
]
\addlegendentry{$j=1$}
\addplot[blue, mark=triangle*, mark options={scale=1}] table[x index=0, y index=1]{
2  0.0976145
4  0.0013239
6  4.27E-05
8  7.13E-07
10  1.1E-08
12  2.32E-10
14  5.23E-12
16  9.21E-14
18  1.19E-15
20  1.12E-15
};
\addlegendentry{$j=3$}
\addplot[green, mark=*, mark options={scale=1}] table[x index=0, y index=1]{
2  0.0903804
4  0.00179832
6  3.07E-05
8  6.39E-07
10  8.02E-09
12  2.03E-10
14  2.96E-12
16  5.69E-14
18  1.43E-15
20  1.24E-15
};
\addlegendentry{$j=5$}
\addplot[red, mark=pentagon*, mark options={scale=1}] table[x index=0, y index=1]{
2  0.0987935
4  0.00219317
6  4.54E-05
8  8.55E-07
10  1.67E-08
12  3.39E-10
14  6.66E-12
16  1.28E-13
18  3.45E-15
20  1.63E-15
};
\addlegendentry{$j=7$}
\addplot[lila, mark=diamond*, mark options={scale=1}] table[x index=0, y index=1]{
2  0.104519
4  0.00247584
6  5.05E-05
8  1.04E-06
10  2.16E-08
12  4.48E-10
14  9.27E-12
16  1.92E-13
18  4.16E-15
20  1.24E-15
};
\end{semilogyaxis}
\end{tikzpicture}
\caption{\label{fig:MassComp}Comparison between 
the convergence rates of the matrix square root with 
respect to the number $\hat{K}$ of terms in \eqref{eq:sqrt}
for different polynomial degrees $p$ and fixed level 
$j = 6$ (left) and different refinement levels $j$ and 
fixed polynomial degree $p=3$ (right).}  
\end{center}
\end{figure}
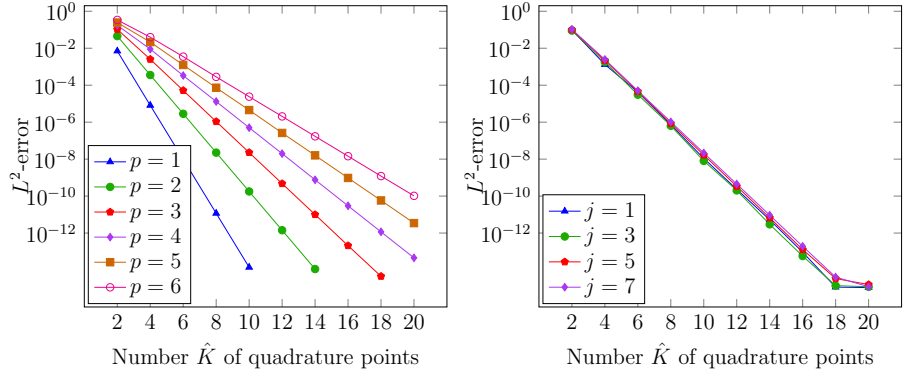

We see in Figure \ref{fig:MassComp} that the 
convergence of the approximation towards the desired 
solution is indeed exponential in the parameter $\hat{K}$, 
with constants depending on the condition number of the 
mass matrix (left plot in Figure \ref{fig:MassComp}), thus 
on the polynomial degree $p$, but not on the refinement 
level $j$ (right plot in Figure \ref{fig:MassComp}). This 
results in linear complexity in the number of degrees of 
freedom for a fixed value of $\hat{K}$ and thus in a complexity 
$\mathcal{O}(N \log N)$ for the desired accuracy matching 
the discretization error of the ansatz space.

\subsubsection*{Gaussian random field simulation}
\begin{figure}[hbt]
  \centering
  \begin{subfigure}[b]{0.33\textwidth}
    \centering
    \includegraphics[width=0.9\linewidth]{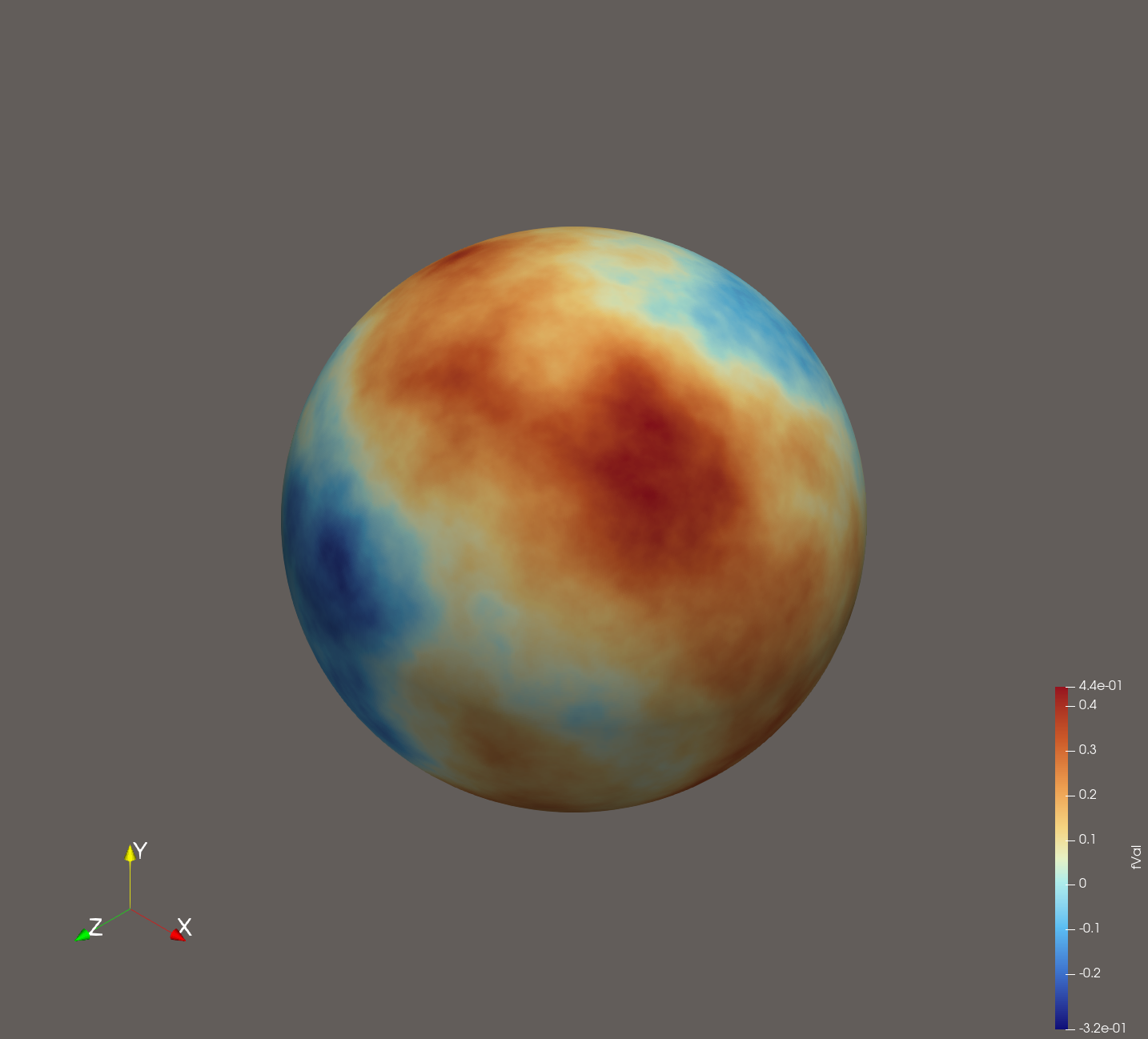}
    \caption{$\beta = 1$ and $\kappa = 1$.}
  \end{subfigure}%
  \begin{subfigure}[b]{0.33\textwidth}
    \centering
    \includegraphics[width=0.9\linewidth]{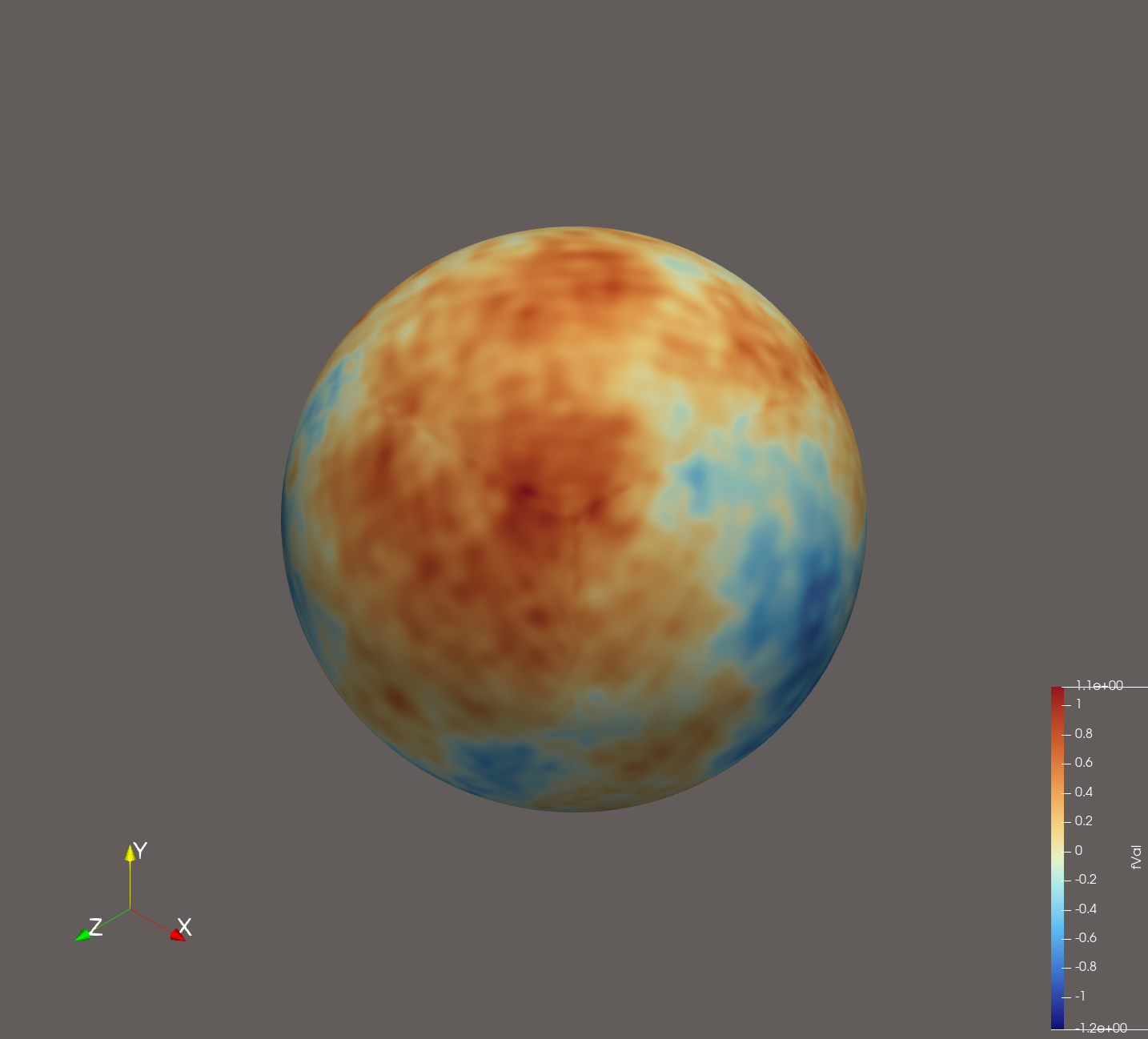}
    \caption{$\beta = 0.75$ and $\kappa = 1$.}
  \end{subfigure}%
  \begin{subfigure}[b]{0.33\textwidth}
    \centering
    \includegraphics[width=0.9\linewidth]{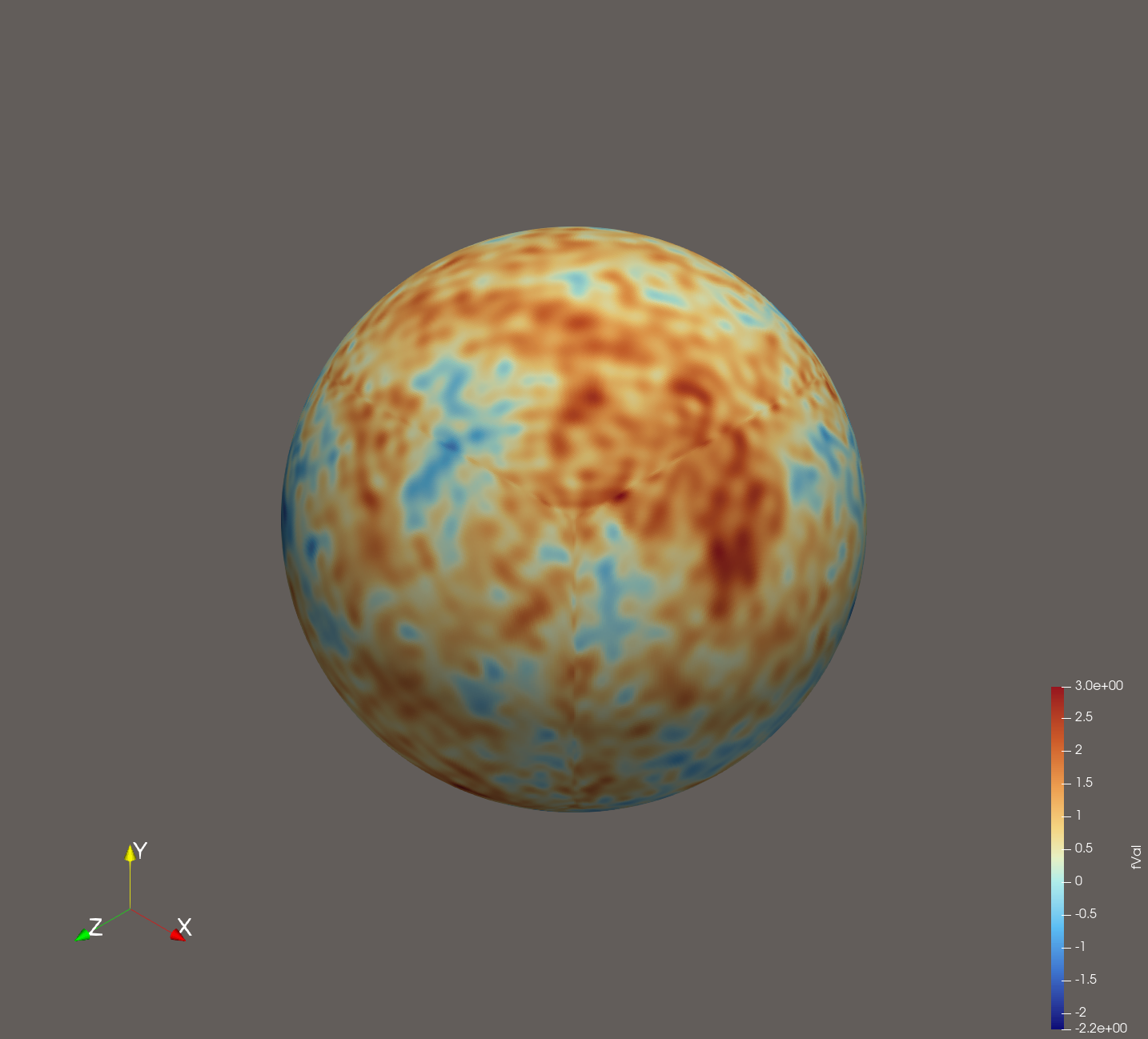}
    \caption{$\beta = 0.55$ and $\kappa = 1$.}
  \end{subfigure}\\
  \begin{subfigure}[b]{0.33\textwidth}
    \centering
    \includegraphics[width=0.9\linewidth]{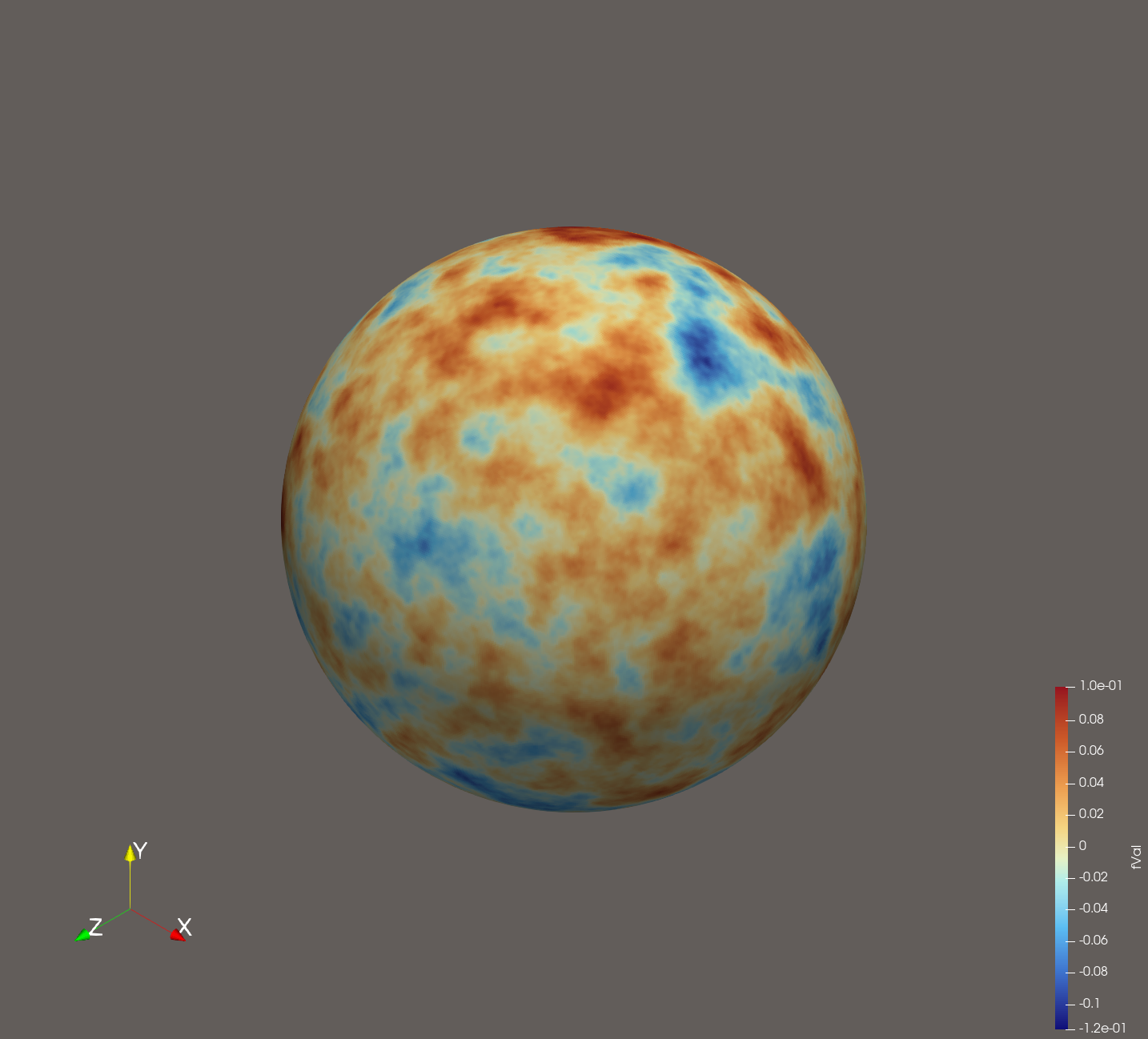}
    \caption{$\beta = 1$ and $\kappa = 10$.}
  \end{subfigure}%
  \begin{subfigure}[b]{0.33\textwidth}
    \centering
    \includegraphics[width=0.9\linewidth]{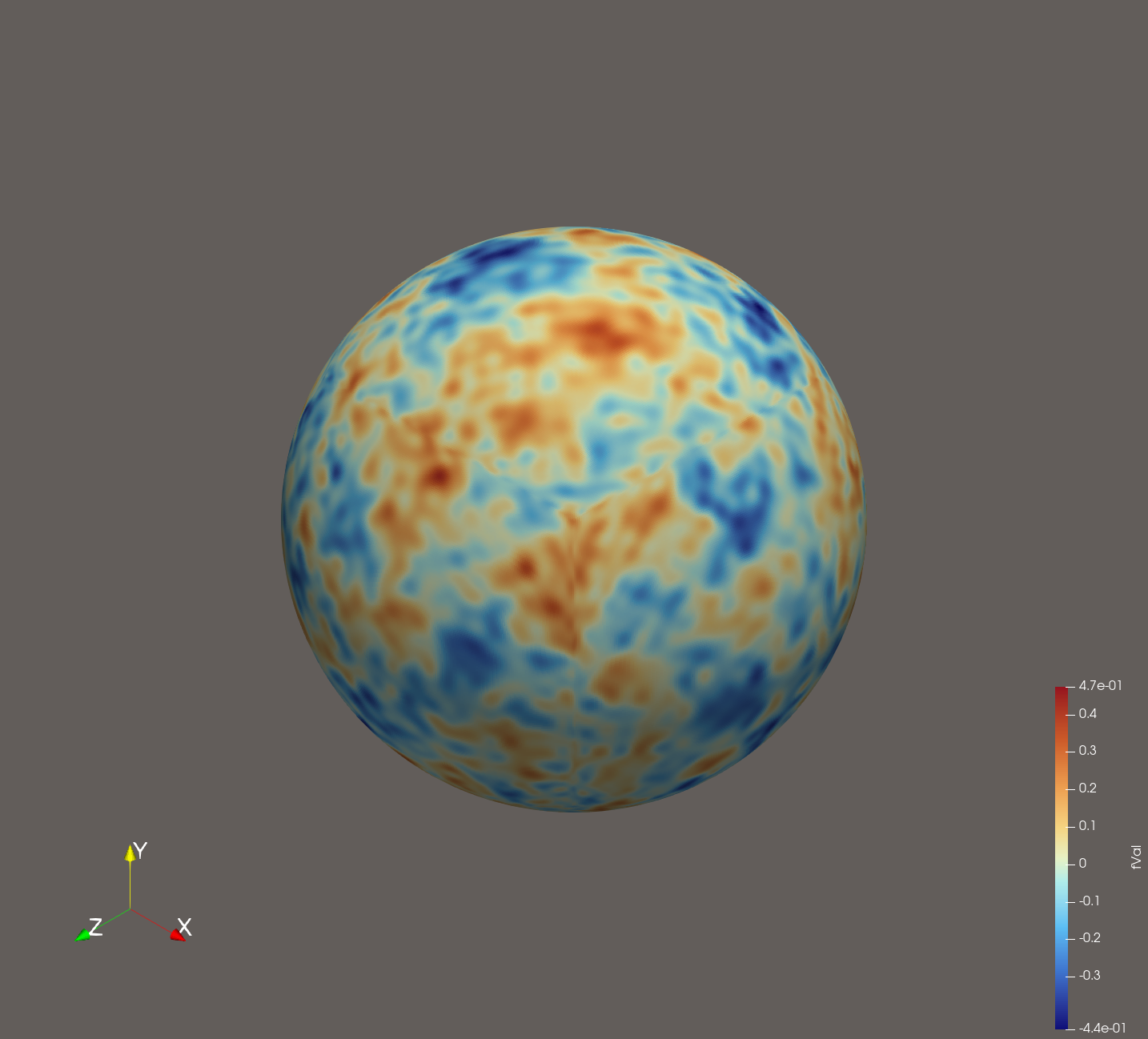}
    \caption{$\beta = 0.75$ and $\kappa = 10$.}
  \end{subfigure}%
  \begin{subfigure}[b]{0.33\textwidth}
    \centering
    \includegraphics[width=0.9\linewidth]{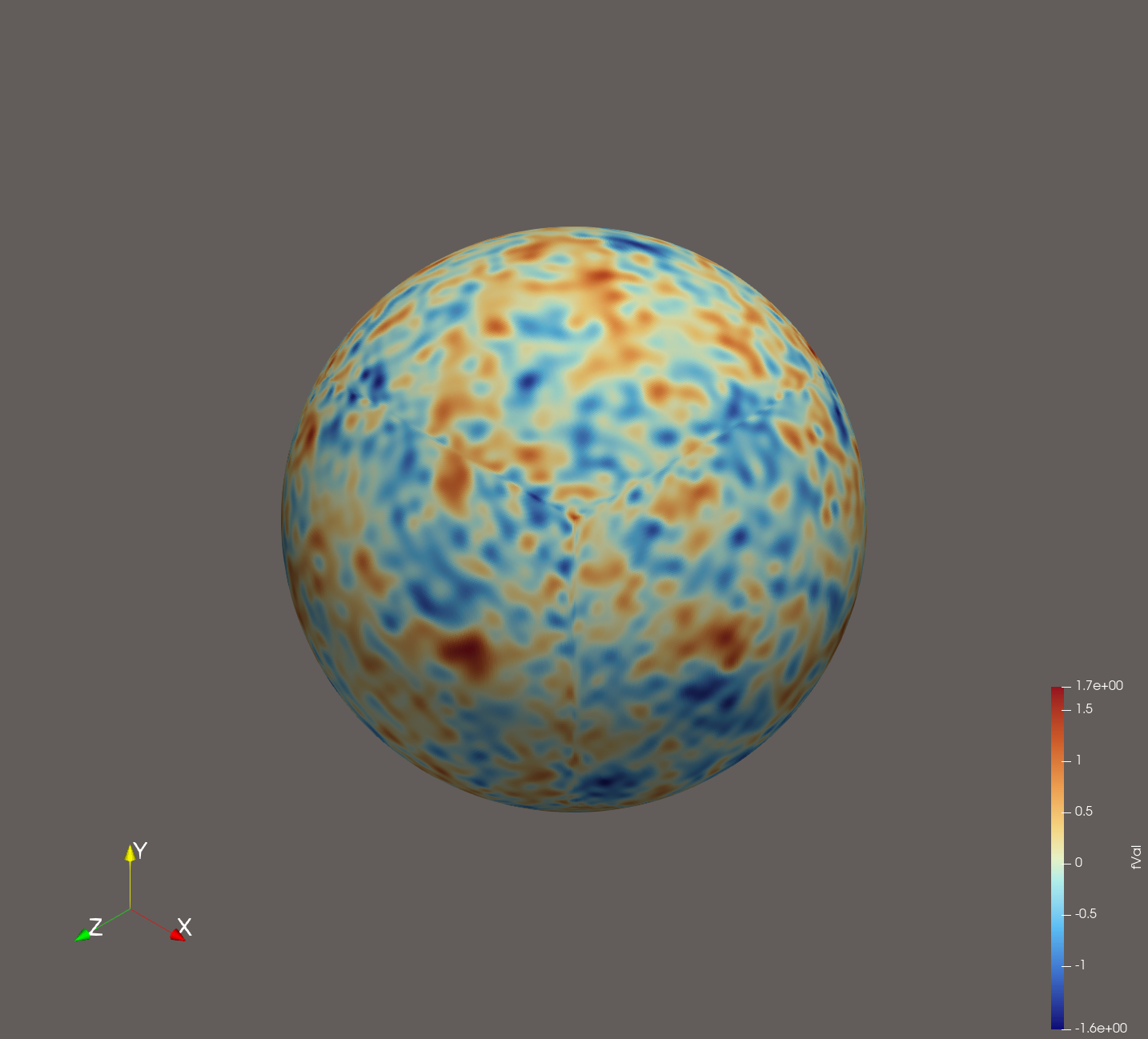}
    \caption{$\beta = 0.55$ and $\kappa = 10$.}
  \end{subfigure}%
  \caption{\label{fig:GRF-sphere}Whittle-Mat\'ern Gaussian random fields on 
  the unit sphere for different parameters of the fractional index $\beta$ and 
  the correlation length $\kappa$.}
\end{figure}

Last but not least, we present some realizations of Gaussian 
Whittle-Mat\'ern fields for different parameters $\beta$ and 
$\kappa$. As can be seen in Figure~\ref{fig:GRF-sphere}, the 
parameter $\beta$ determines the smoothness of the random 
field while the parameter $\kappa$ is inversely proportional 
to the correlation length of the random field.

\subsection{Drilled cube and Stanford bunny}
\begin{figure}[hbt]
  \centering
  \includegraphics[width=0.45\textwidth]{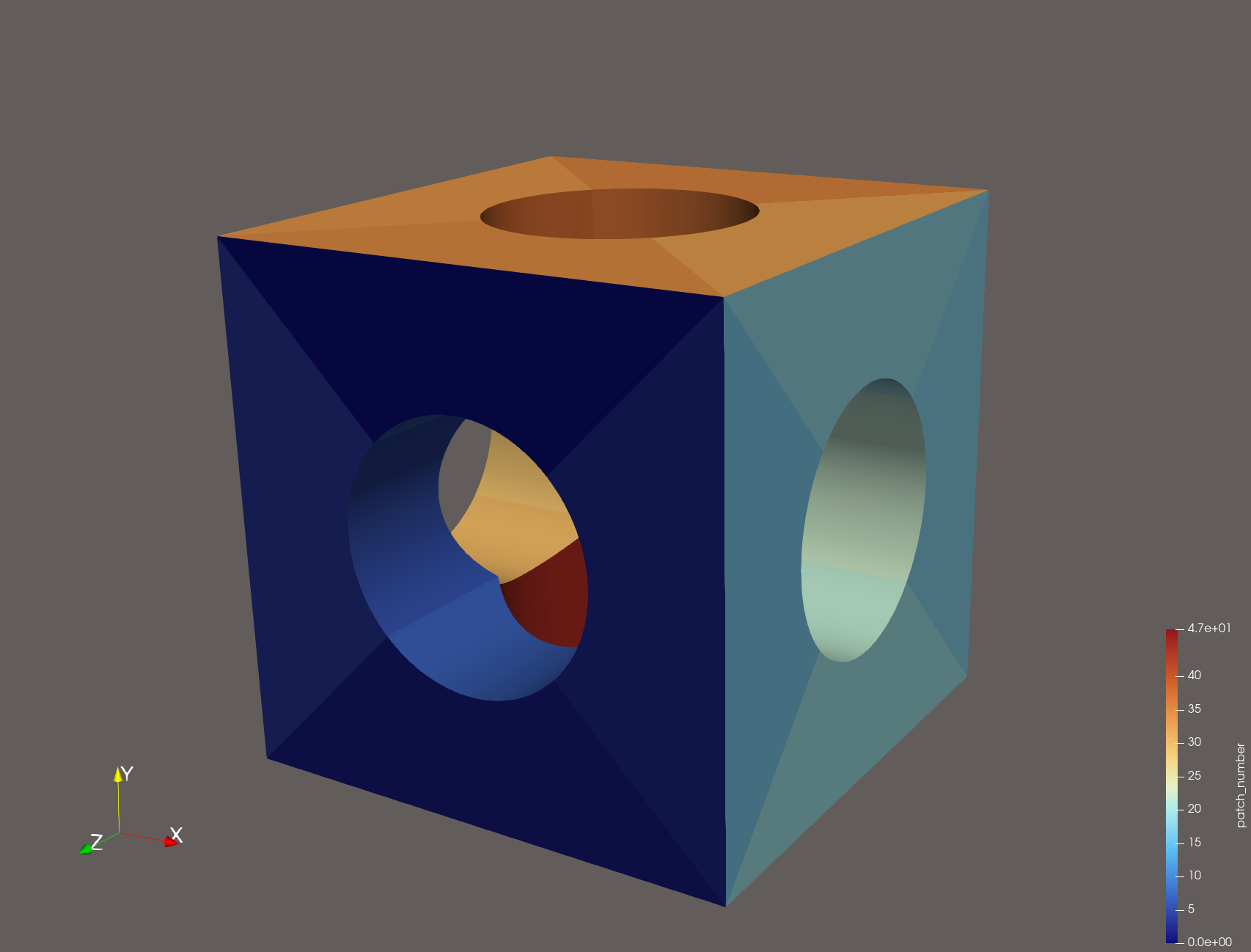}
  \includegraphics[width=0.45\textwidth,trim={0 0 0 205},clip]{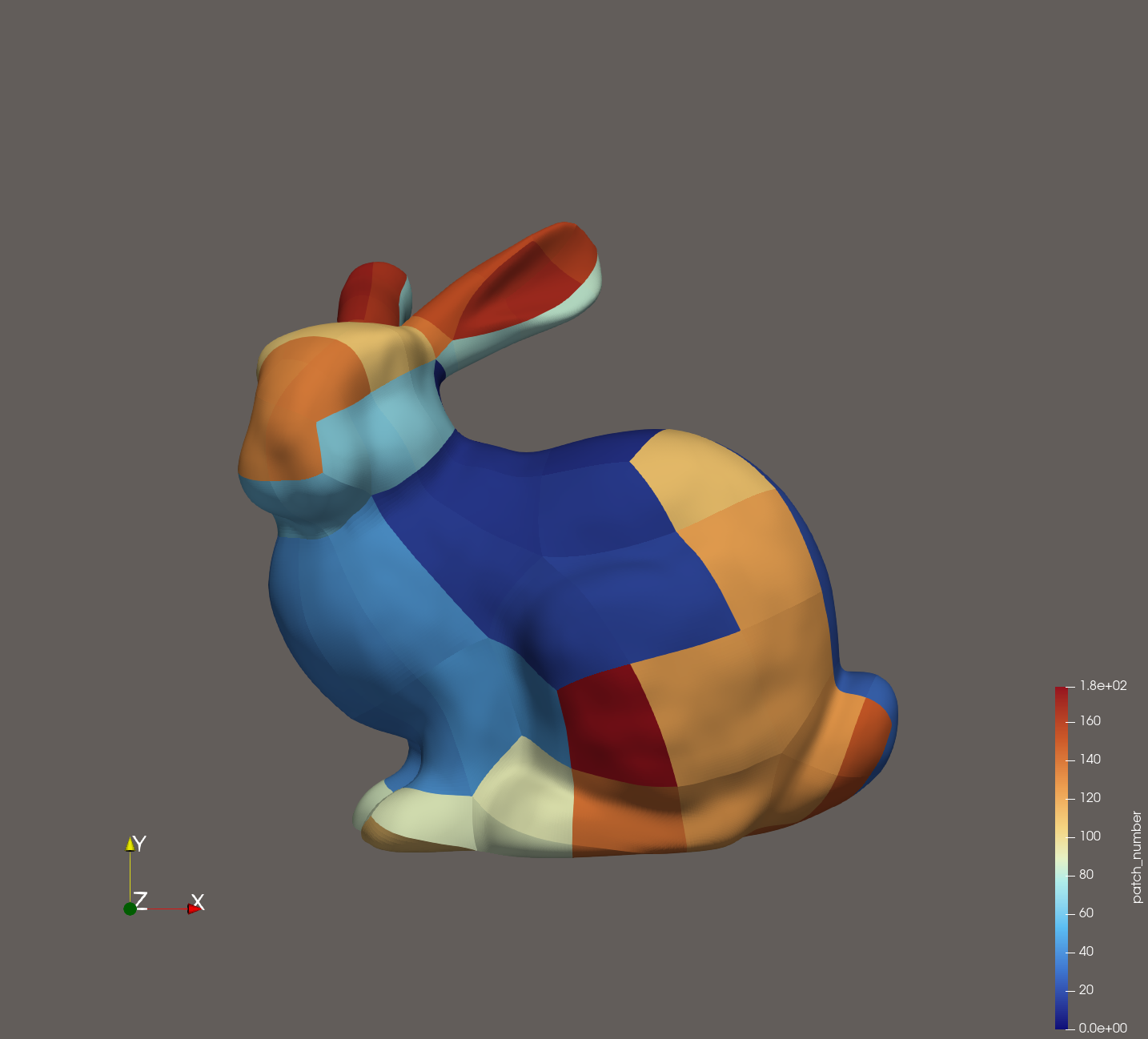}
  \caption{\label{fig:parametrizations}Parametrizations of a drilled cube (left) 
  by 48 quadrangular patches and of the Stanford bunny (right) by 179 
  quadrangular patches.}
\end{figure}

Finally, we demonstrate that our approach applies also 
to other, especially nontrivial surface representations: we turn 
our attention to a drilled cube and to the famous Stanford bunny found 
in Figure~\ref{fig:parametrizations}. The drilled cube is parametrized by
48 quadrangular patches while the Stanford bunny is parametrized 
by 179 patches. The parametrization of the Stanford bunny has 
been introduced in \cite{huang_isogeometric_2022} and is based 
on the quadrangulation from \cite{pietroni_tracing_2016}. We 
refer to these references for further details about the generation
of the surface representation.

Although the surface of the Stanford bunny is much more complex
in comparison with the unit sphere, all the algorithms behave essentially 
the same as on the unit sphere, such that we refrain from presenting 
detailed results here. In particular, the performance of the V-cycle 
multigrid is rather poor in case of the Stanford bunny since we need 
a lot of pre- and post-smoothing steps to ensure convergence. 
In contrast, the preconditioned conjugate gradient method 
converges well for all polynomial degrees $1\le p \le 5$ in case 
of both BPX-variants. Table~\ref{tab:bunny_BPX} reports the 
results for the Stanford bunny, which are quite similar to that
obtained for the unit sphere.

\begin{table}[hbt]
\begin{center}
{\small
\begin{tabular}{|c|ccccc|ccccc|}\cline{2-11}
\multicolumn{1}{c|}{}
& \multicolumn{5}{c|}{BPX with diagonal scaling} & \multicolumn{5}{c|}{BPX with SSOR}\\\hline
$j$ & $p=1$ & $p=2$ & $p=3$ & $p=4$ & $p=5$ &$p=1$ & $p=2$ & $p=3$ & $p=4$ & $p=5$ \\\hline
0 & 64 & 89 & 143 & 261 & 531 & 27 & 35 & 47 & 73 & 150 \\
1 & 79 & 123 & 190 & 335 & 751 & 40 & 49 & 61 & 95 & 183 \\
2 & 100 & 144 & 216 & 347 & 754 & 53 & 60 & 72 & 102 & 198 \\
3 & 132 & 170 & 244 & 376 & 690 & 67 & 69 & 85 & 116 & 178 \\
4 & 178 & 203 & 278 & 425 & 708 & 85 & 82 & 98 & 138 & 198 \\
5 & 230 & 236 & 322 & 485 & 814 & 104 & 94 & 115 & 163 & 245 \\\hline
\end{tabular}}
\caption{\label{tab:bunny_BPX}Number of iterations of the preconditioned 
conjugate gradient with BPX-preconditioning in case of a diagonal scaling 
on each level and in case of an SSOR preconditioning on each level.}
\end{center}
\end{table}

\begin{figure}[hbt]
  \centering
  \begin{subfigure}[b]{0.33\textwidth}
    \centering
    \includegraphics[width=0.9\linewidth]{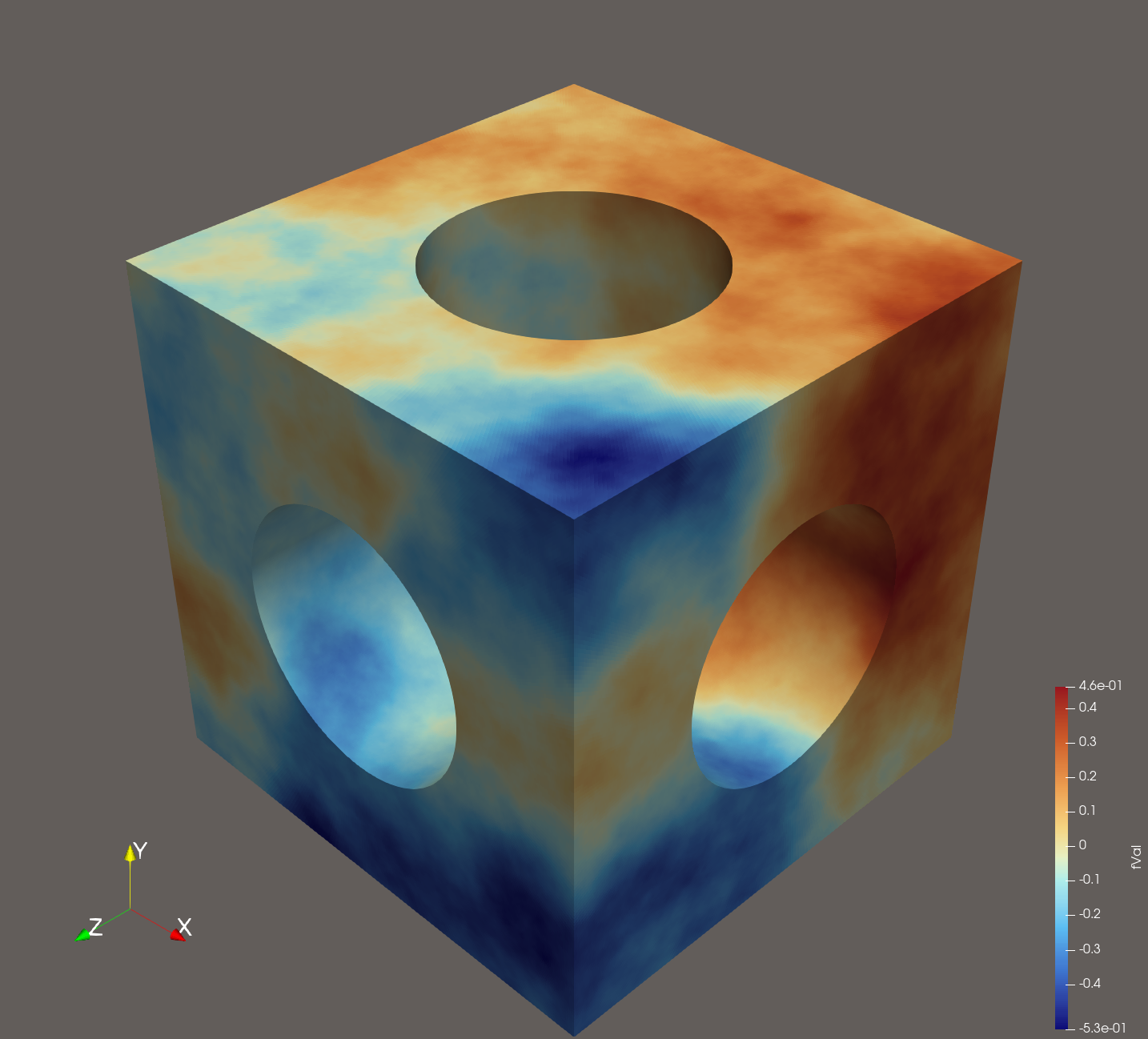}
    \caption{$\beta = 1$ and $\kappa = 1$.}
  \end{subfigure}%
  \begin{subfigure}[b]{0.33\textwidth}
    \centering
    \includegraphics[width=0.9\linewidth]{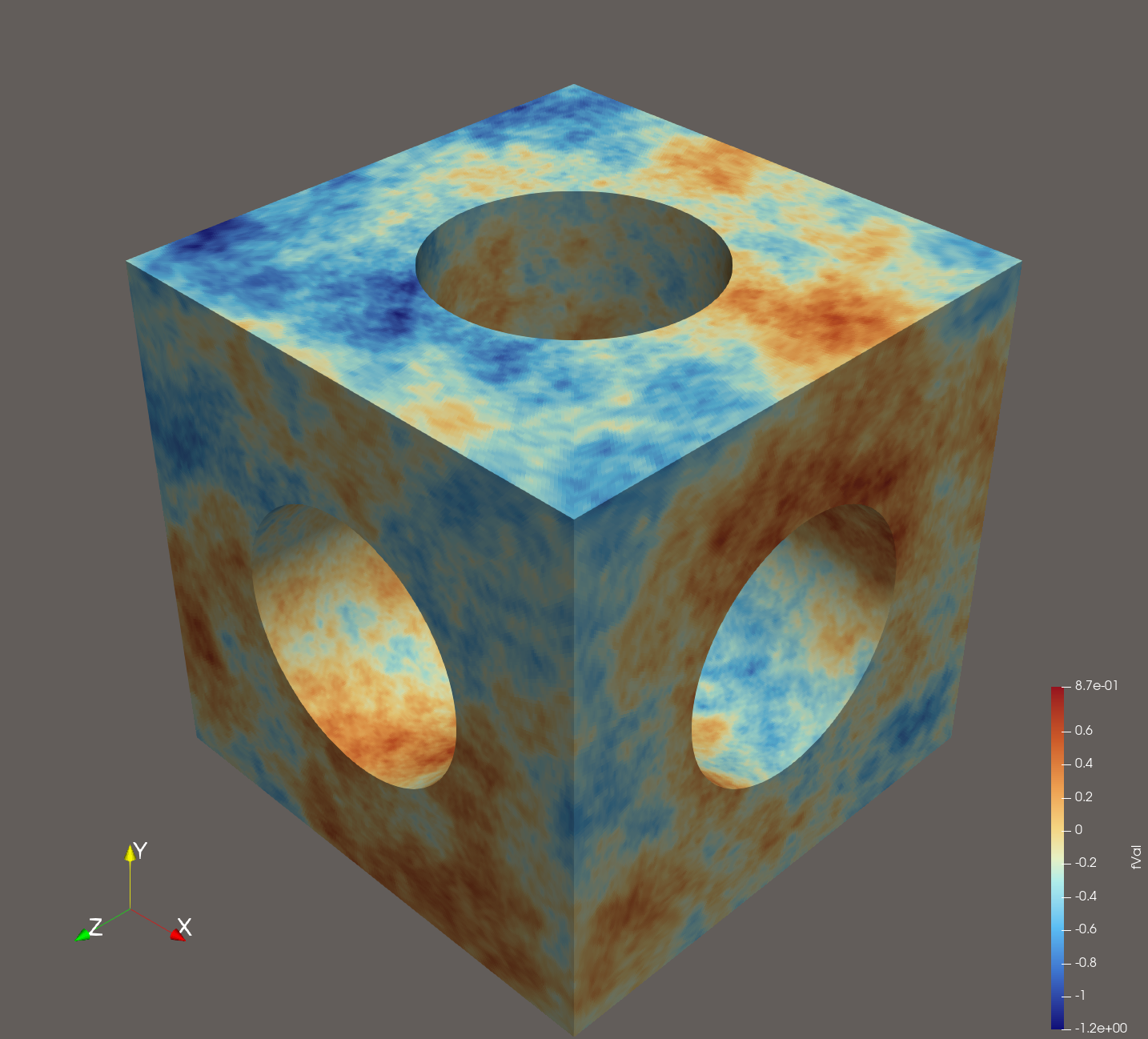}
    \caption{$\beta = 0.75$ and $\kappa = 1$.}
  \end{subfigure}%
  \begin{subfigure}[b]{0.33\textwidth}
    \centering
    \includegraphics[width=0.9\linewidth]{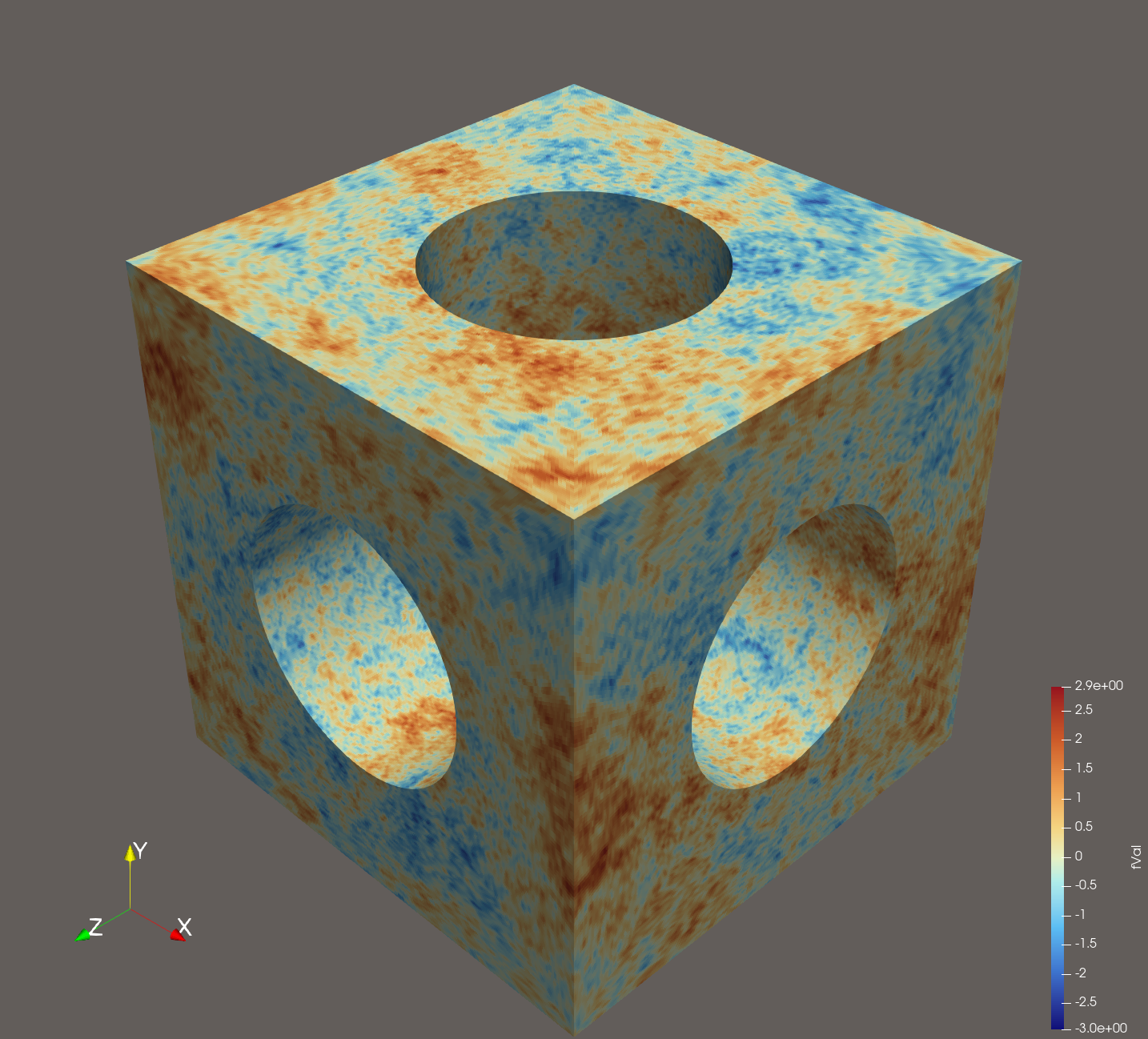}
    \caption{$\beta = 0.55$ and $\kappa = 1$.}
  \end{subfigure}\\
  \begin{subfigure}[b]{0.33\textwidth}
    \centering
    \includegraphics[width=0.9\linewidth]{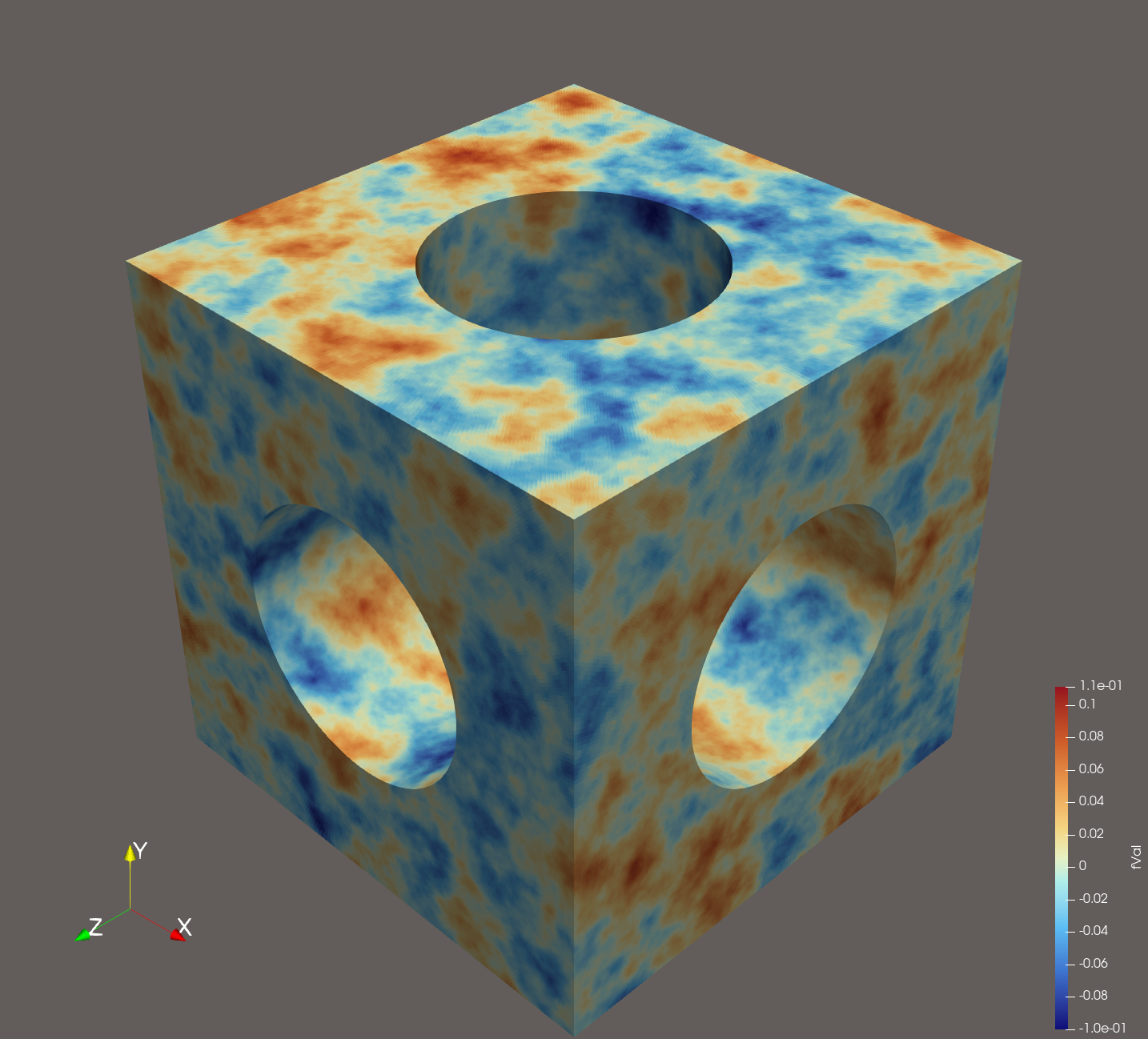}
    \caption{$\beta = 1$ and $\kappa = 10$.}
  \end{subfigure}%
  \begin{subfigure}[b]{0.33\textwidth}
    \centering
    \includegraphics[width=0.9\linewidth]{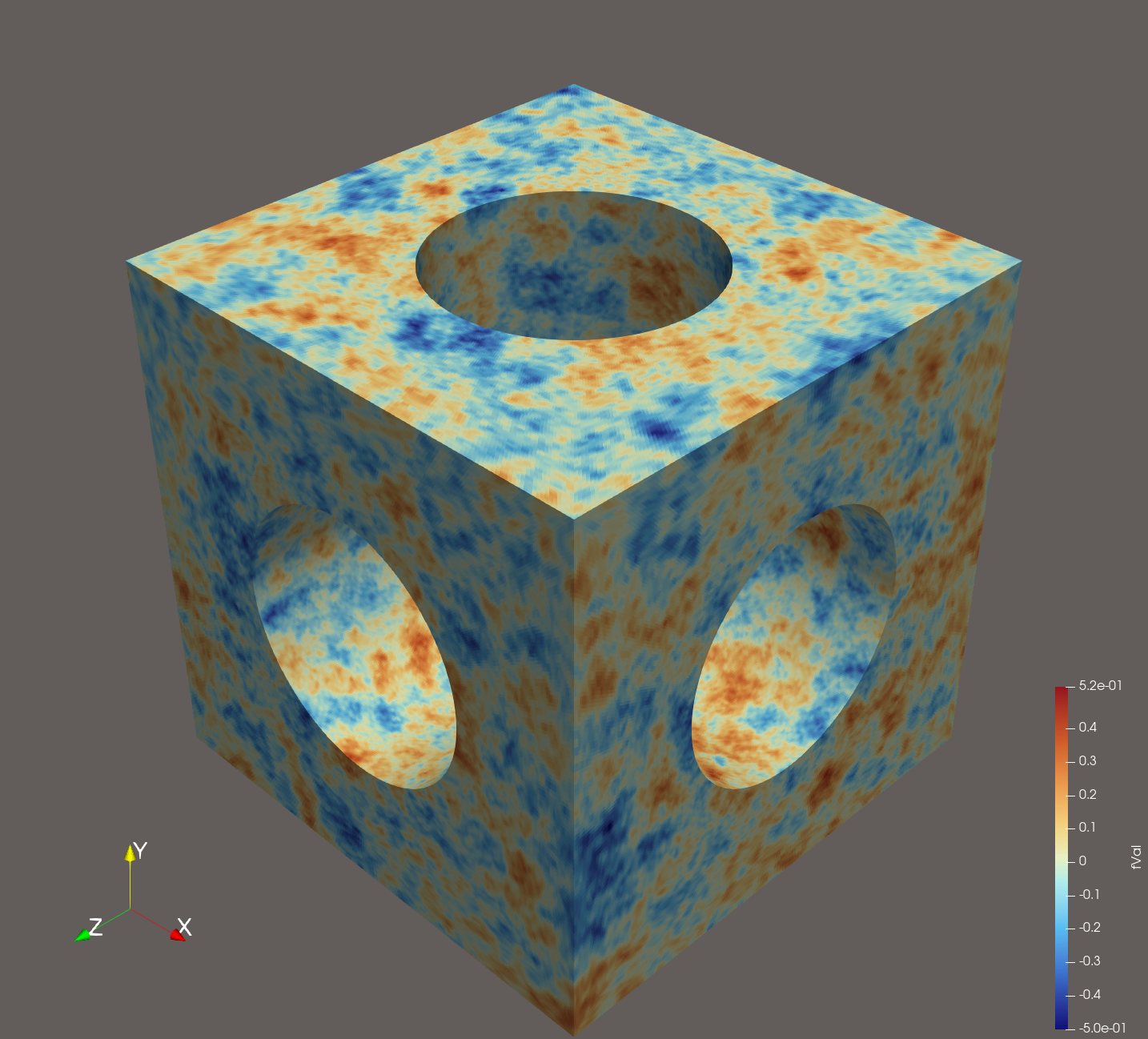}
    \caption{$\beta = 0.75$ and $\kappa = 10$.}
  \end{subfigure}%
  \begin{subfigure}[b]{0.33\textwidth}
    \centering
    \includegraphics[width=0.9\linewidth]{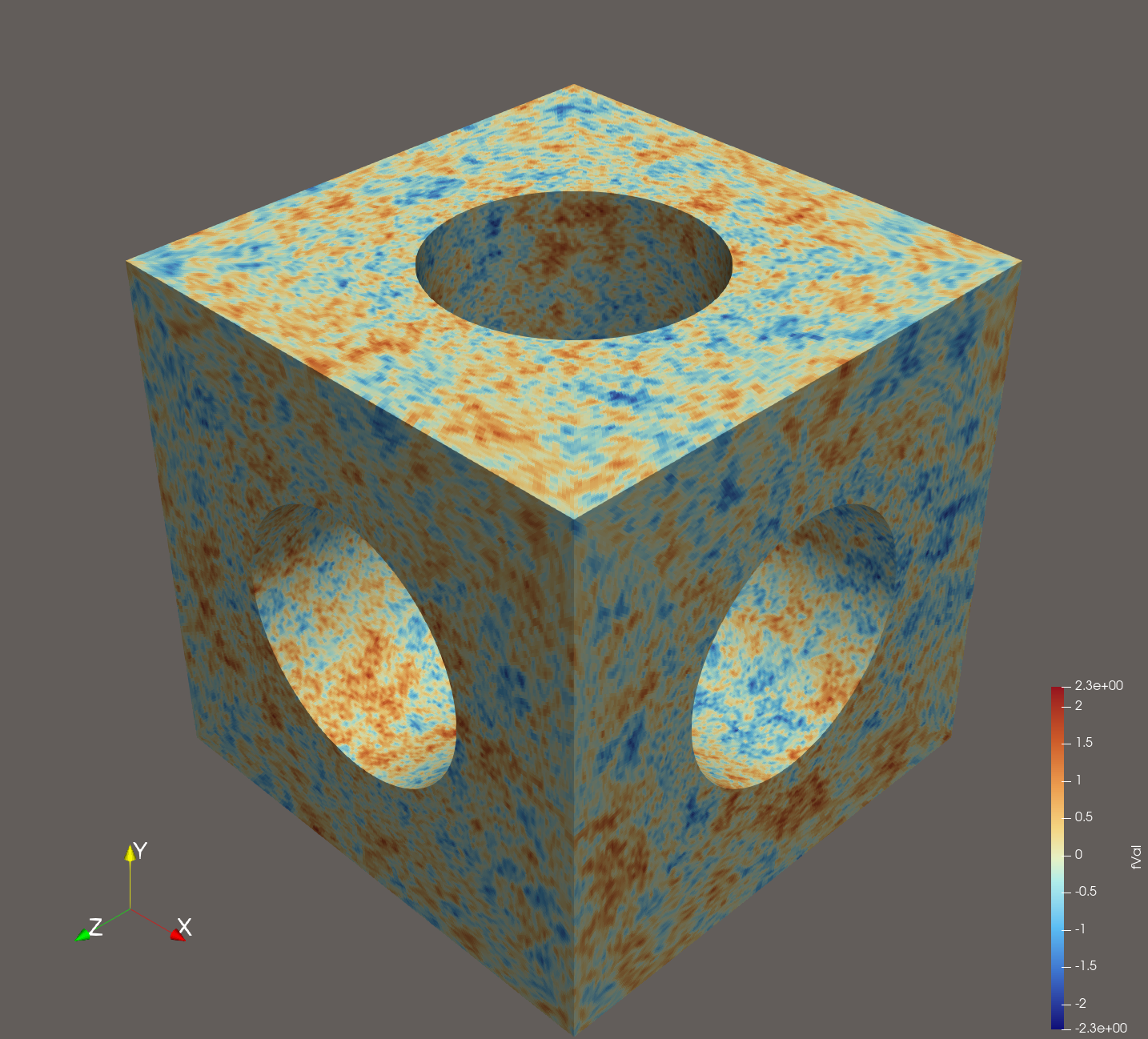}
    \caption{$\beta = 0.55$ and $\kappa = 10$.}
  \end{subfigure}%
  \caption{\label{fig:cube}Whittle-Mat\'ern Gaussian random fields on 
  drilled cube for different parameters of the fractional index $\beta$ and 
  the correlation length $\kappa$.}
\end{figure}

\begin{figure}[hbt]
  \centering
  \begin{subfigure}[b]{0.33\textwidth}
    \centering
    \includegraphics[width=0.9\linewidth,trim={0 0 0 205},clip]{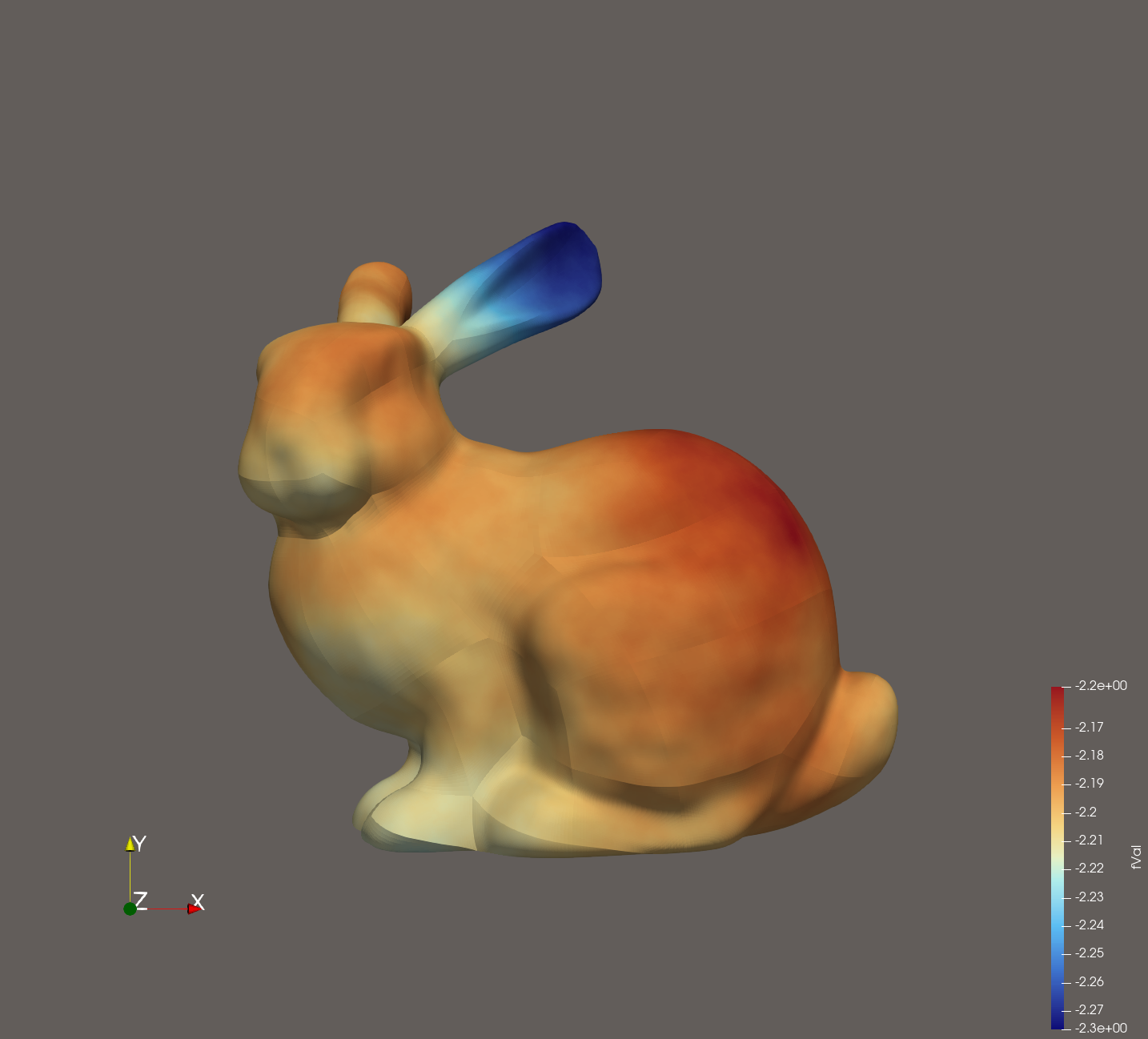}
    \caption{$\beta = 1$ and $\kappa = 1$.}
  \end{subfigure}%
  \begin{subfigure}[b]{0.33\textwidth}
    \centering
    \includegraphics[width=0.9\linewidth,trim={0 0 0 205},clip]{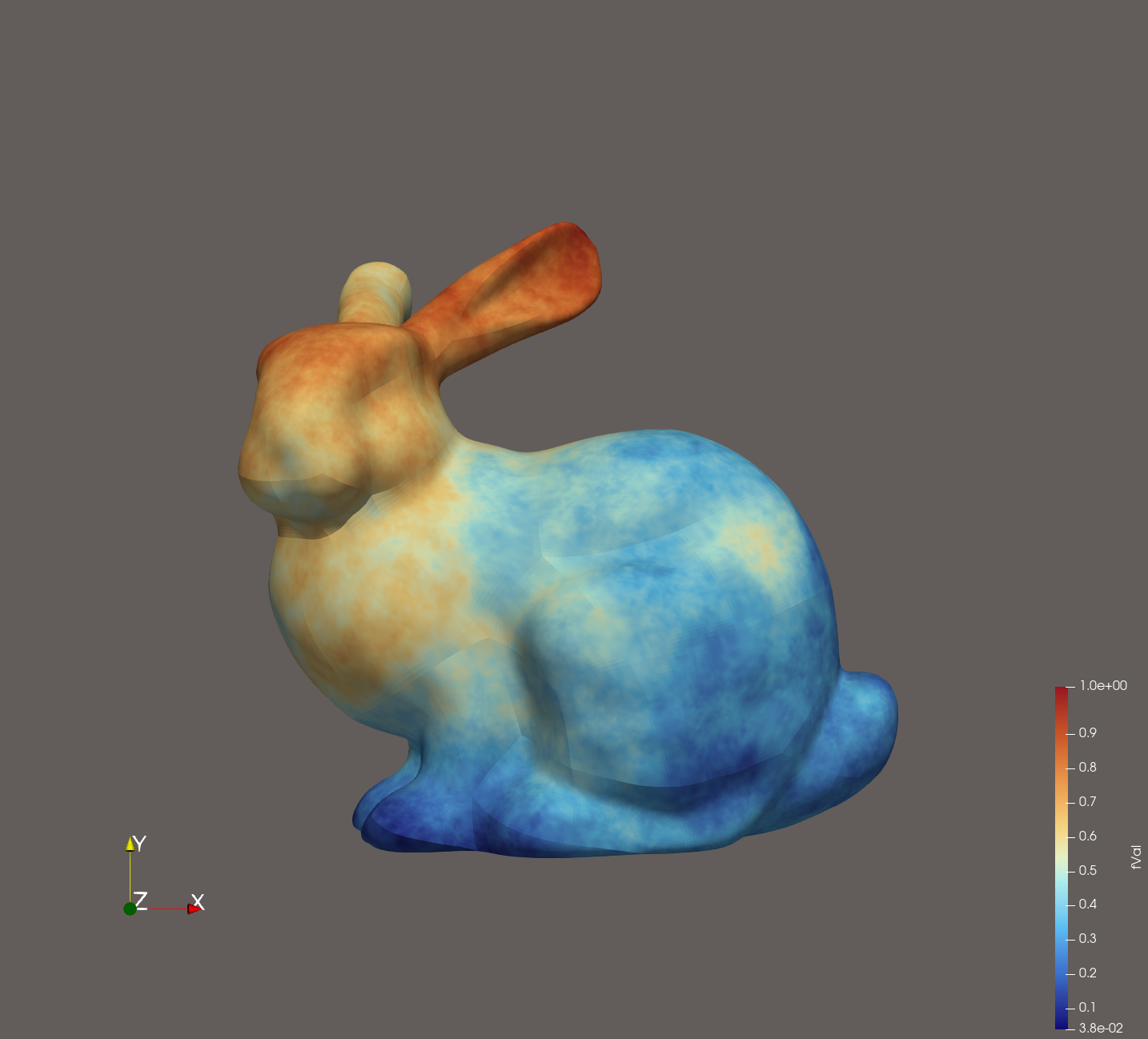}
    \caption{$\beta = 0.75$ and $\kappa = 1$.}
  \end{subfigure}%
  \begin{subfigure}[b]{0.33\textwidth}
    \centering
    \includegraphics[width=0.9\linewidth,trim={0 0 0 205},clip]{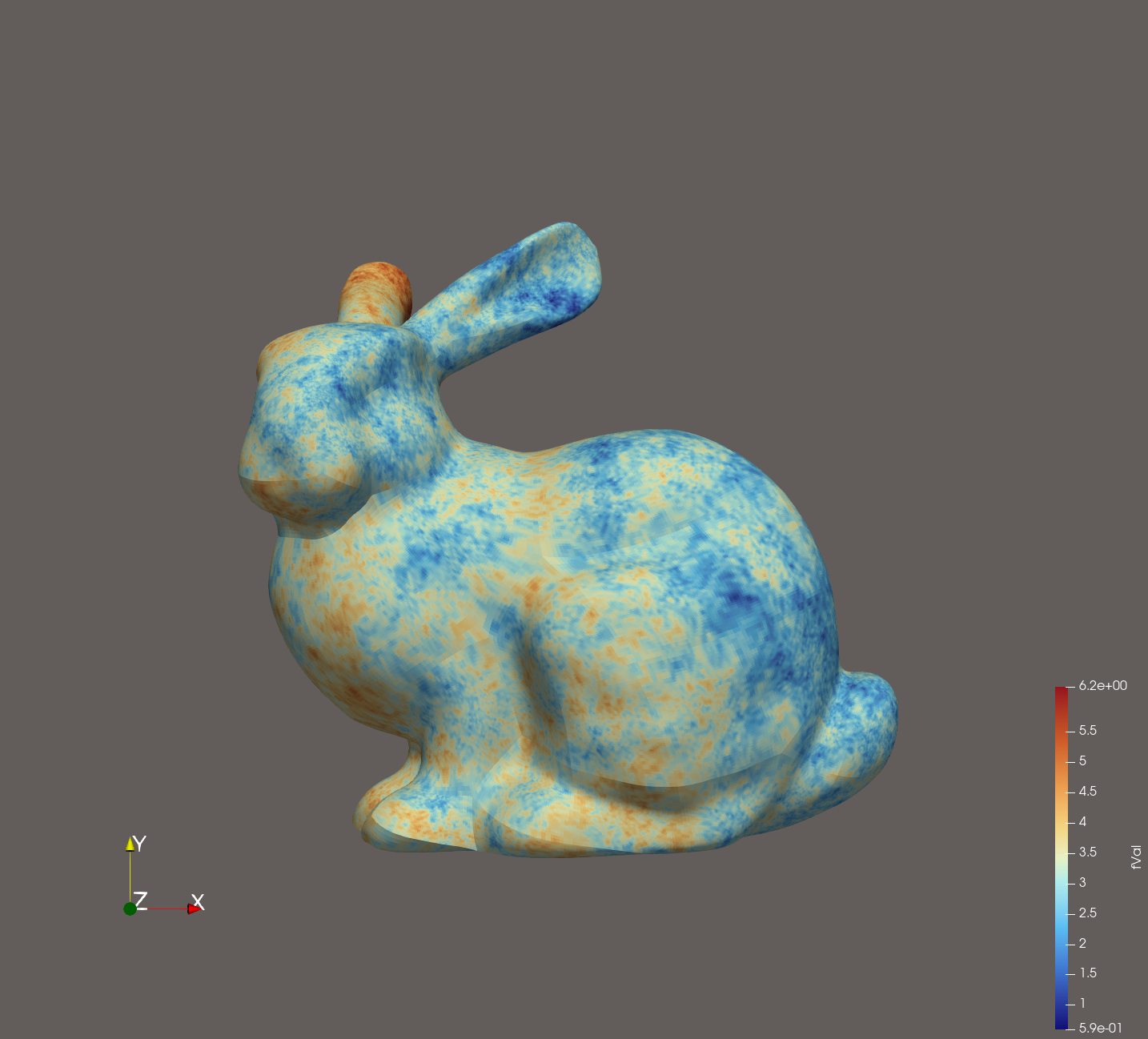}
    \caption{$\beta = 0.55$ and $\kappa = 1$.}
  \end{subfigure}\\
  \begin{subfigure}[b]{0.33\textwidth}
    \centering
    \includegraphics[width=0.9\linewidth,trim={0 0 0 205},clip]{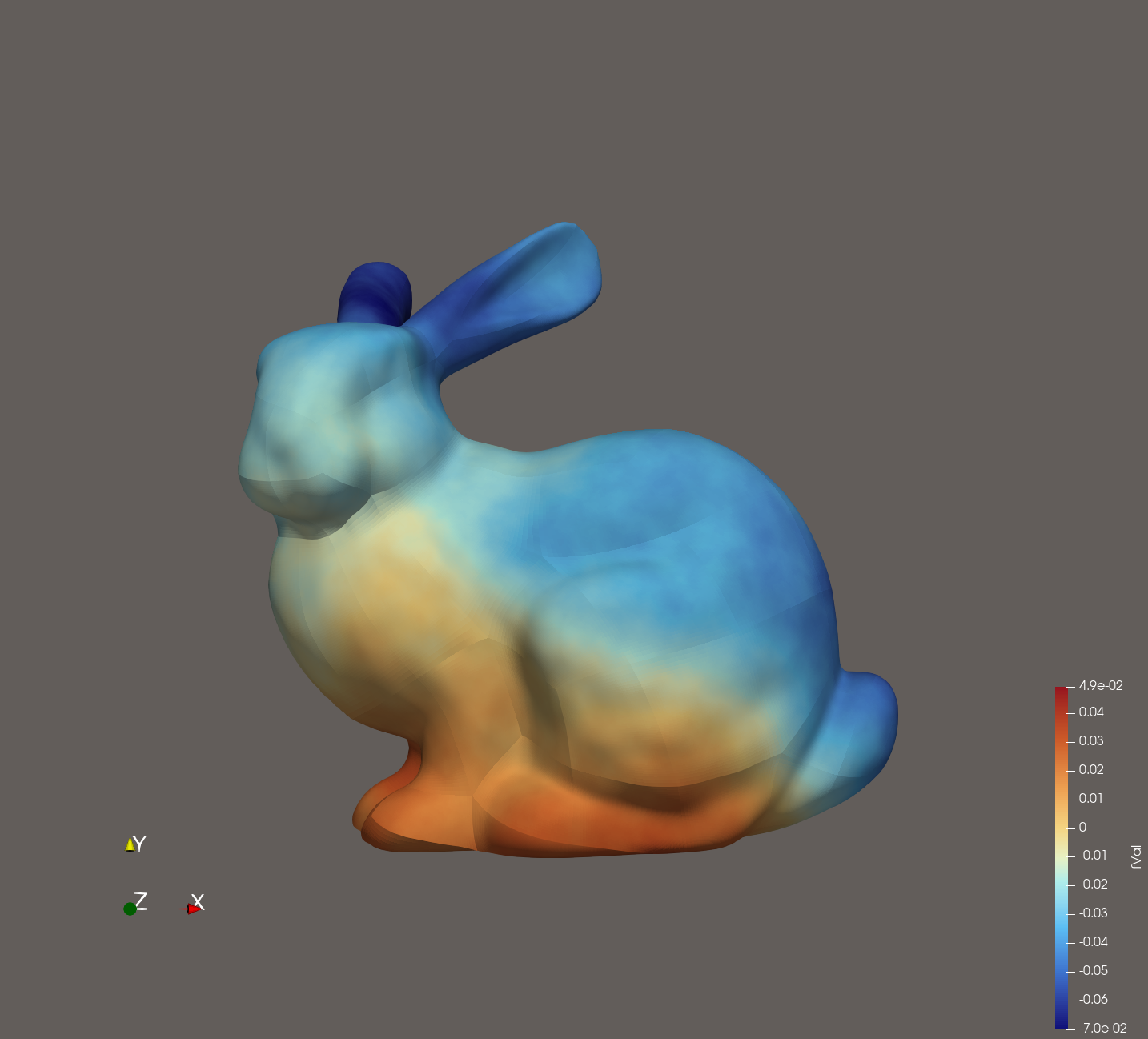}
    \caption{$\beta = 1$ and $\kappa = 10$.}
  \end{subfigure}%
  \begin{subfigure}[b]{0.33\textwidth}
    \centering
    \includegraphics[width=0.9\linewidth,trim={0 0 0 205},clip]{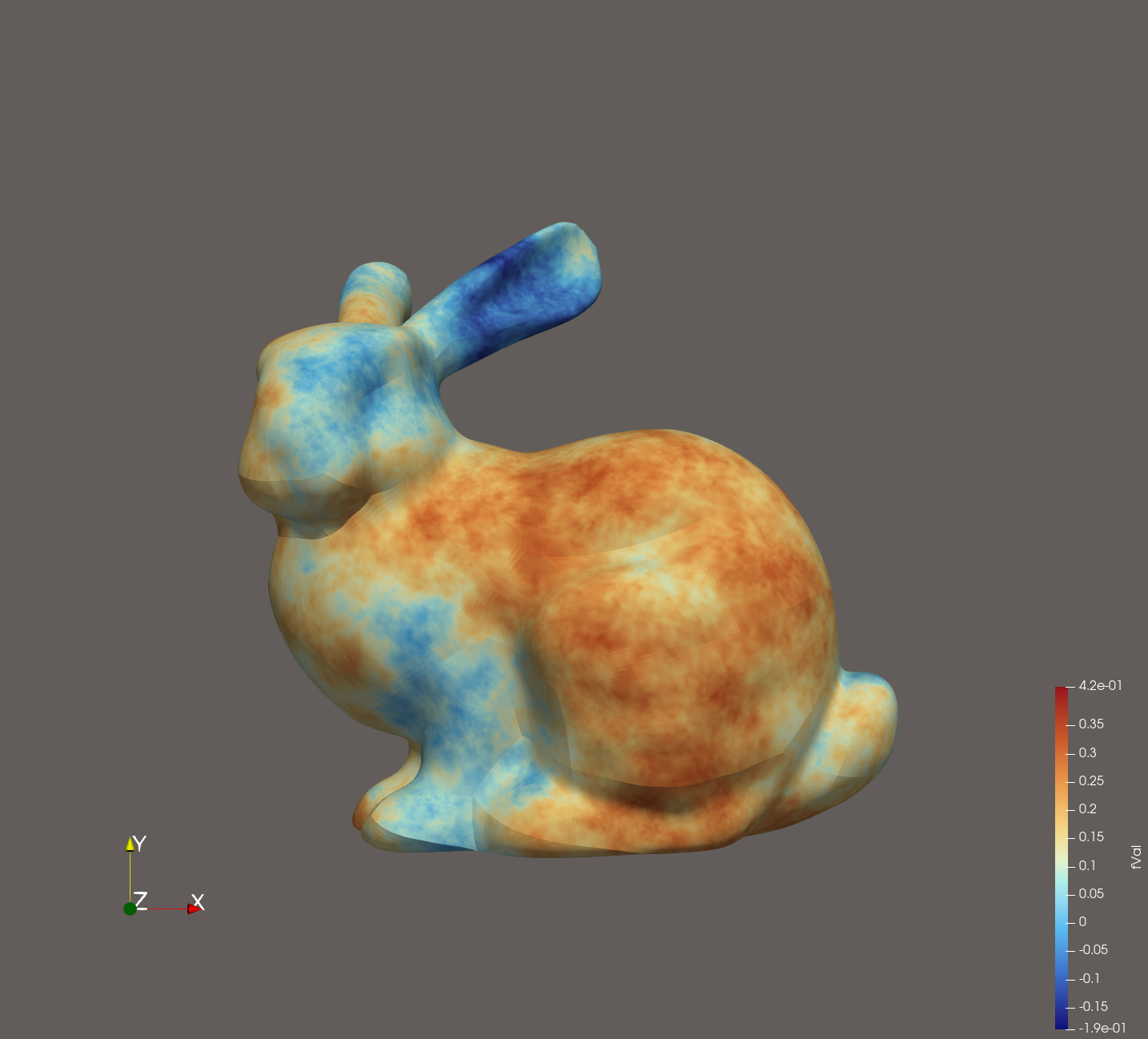}
    \caption{$\beta = 0.75$ and $\kappa = 10$.}
  \end{subfigure}%
  \begin{subfigure}[b]{0.33\textwidth}
    \centering
    \includegraphics[width=0.9\linewidth,trim={0 0 0 205},clip]{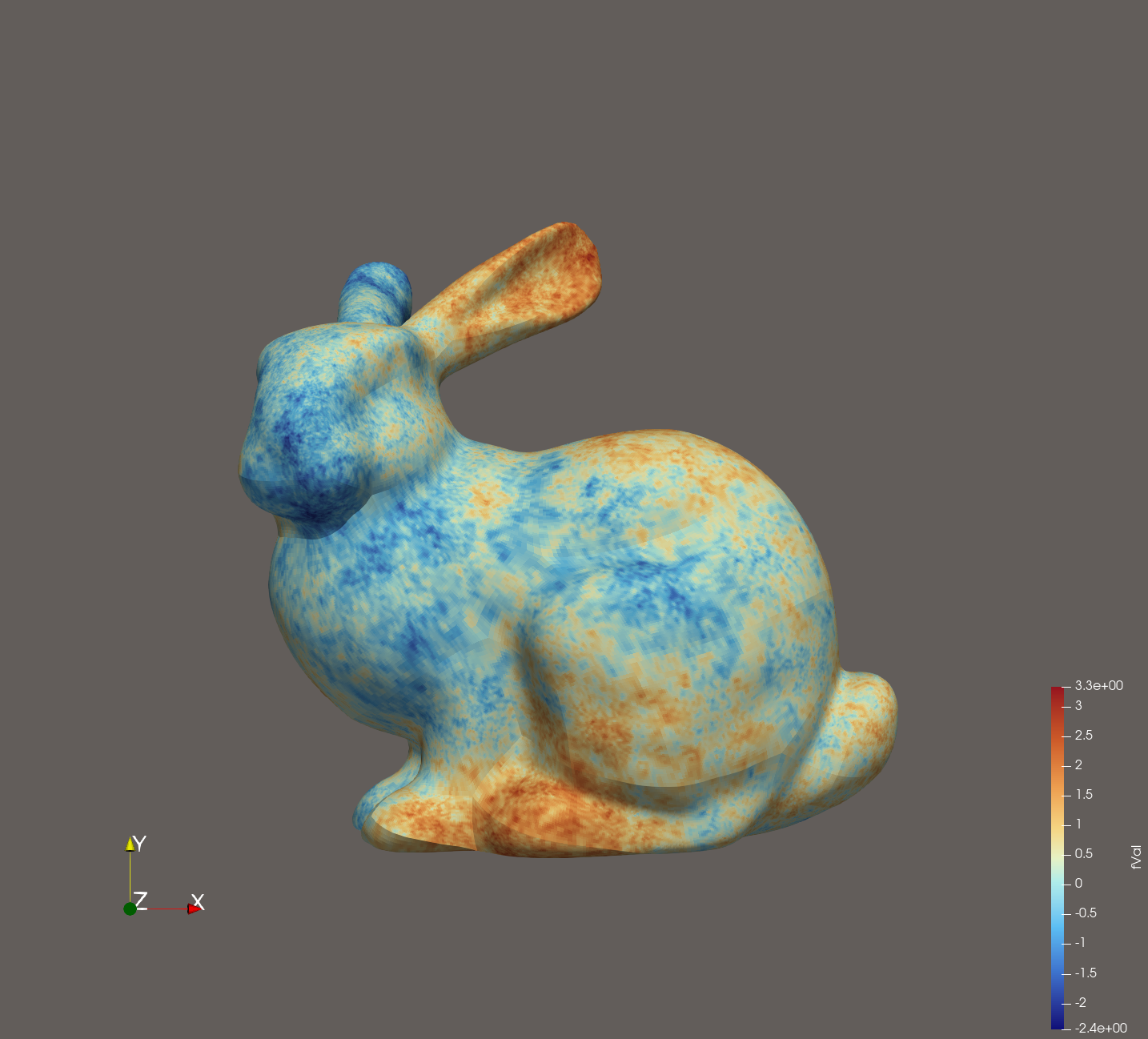}
    \caption{$\beta = 0.55$ and $\kappa = 10$.}
  \end{subfigure}%
  \caption{\label{fig:bunny}Whittle-Mat\'ern Gaussian random fields on 
  the Stanford bunny for different parameters of the fractional index 
  $\beta$ and the correlation length $\kappa$.}
\end{figure}

On both surfaces under consideration, we simulate 
Gaussian random fields generated by Whittle-Mat\'ern 
covariance functions with the same set of parameters as 
for the unit sphere. Visualizations of these Whittle-Mat\'ern 
Gaussian random fields can be found in Figure~\ref{fig:cube}
for the drilled cube and in Figure~\ref{fig:bunny} for the Stanford 
bunny. The level of refinement has been set to $j=5$ and 
the polynomial degree to $p = 2$.

\section{Conclusion}
\label{sct:conclusio}
In the present article, we proposed the application of 
the isogeometric finite element method for the fast simulation
of Gaussian random fields on surfaces which are generated 
by the (Whittle-) Mat\'ern class of covariance functions. The
isogeometric finite element method is able to compute realizations 
of the desired Gaussian random field in essentially linear overall 
complexity. This means the complexity is linear in the number of 
degrees of freedom except for polylogarithmic terms. 

Compared to previous approaches from \cite{Bonito,Annika2,
Annika1} which are based on planar or curved \emph{linear} 
surface finite elements, the isogeometric finite element method 
enables higher rates of convergences in case of Gaussian 
random fields of sufficiently high smoothness since the 
polynomial degree $p$ of the ansatz spaces can also be 
chosen larger than $p=1$. Numerical results have been 
reported which validate the applicability and feasibility 
of the isogeometric approach also on complicated 
surfaces. 

Although we restricted ourselves to stationary Gaussian 
random fields being defined on two-dimensional surfaces, 
the extension of the present methodology to nonstationary 
Gaussian random fields is easily possible, see 
e.g.~\cite{Annika1,Rue1} for the details.

\bibliographystyle{plain}

\end{document}